\documentclass[12pt]{amsart}
\usepackage{amsmath,amsfonts,amssymb}

\usepackage{stmaryrd}
\usepackage{bbding}
\usepackage{amsfonts}
\usepackage[colorlinks,linkcolor=blue,urlcolor=cyan,citecolor=green,pagebackref]{hyperref}
\usepackage{amsmath,amssymb,amscd,amsthm}
\allowdisplaybreaks
\usepackage{pstricks,pst-node}

\usepackage{mathrsfs}
\usepackage{url}
\usepackage{amssymb}
\usepackage{latexsym}
\usepackage{amsfonts}
\usepackage{amsmath}
\usepackage{eucal}
\usepackage{bm}
\usepackage{bbm}
\usepackage{graphicx}
\usepackage[english]{varioref}
\usepackage[nice]{nicefrac}
\usepackage[all]{xy}
\usepackage{amsthm}

\newtheorem{Lem}{Lemma}[section]
\newtheorem{Thm}[Lem]{Theorem}
\newtheorem{Cor}{Corollary}[section]
\newtheorem{Rem}{Remark}[section]

\newcommand{\eq}{\begin{equation}}
\newcommand{\en}{\end{equation}}
\newcommand{\beqna}[1]{\begin{eqnarray}\label{#1}}
\newcommand{\eeqna}{\end{eqnarray}}
\newcommand{\beqn}[1]{\begin{equation}\label{#1}}
\newcommand{\eeqn}{\end{equation}}
\numberwithin{equation}{section}
\numberwithin{table}{section}

\begin{document}
\title[]{A characterization of the unitary highest weight modules by Euclidean Jordan algebras}

\author{Zhanqiang Bai}
\address{Department of Mathematics, Hong Kong Univ. of Sci. and Tech., Clear Water Bay, Kowloon, Hong Kong}
\email{mazqbai@ust.hk}

\maketitle

\begin{abstract}
Let  $\mathfrak{co}(J)$ be the  conformal algebra of a simple Euclidean Jordan algebra $J$.  We show that a (non-trivial) unitary highest weight $\mathfrak{co}(J)$-module  has the smallest positive Gelfand-Kirillov dimension
if and only if a certain quadratic relation is satisfied in the universal enveloping algebra  $U(\mathfrak{co}(J)_{\mathbb{C}})$. In particular, we find an quadratic element in $U(\mathfrak{co}(J)_{\mathbb{C}})$. A prime ideal in $U(\mathfrak{co}(J)_{\mathbb{C}})$ equals   the Joseph ideal if and only if it contains this quadratic element.

\end{abstract}
%
%
%

\section{ Introduction}
\subsection{Unitary highest weight modules}
Let ($G$, $K$) be an irreducible Hermitian symmetric pair of non-compact type. If we denote their Lie algebras by ($\mathfrak{g}_{0}$, $\mathfrak{k}_{0}$), denote their complexification by $\mathfrak{g}={\mathfrak{g}_{0}}\bigotimes \mathbb{C}$ and $\mathfrak{k}=\mathfrak{k}_{0}\bigotimes\mathbb{C}$, then $\mathfrak{k}=\mathbb{C}H\bigoplus[\mathfrak{k}, \mathfrak{k}]$ with $ad(H)$ having eigenvalues $0, 1, -1$ on $\mathfrak{g}$. Usually we use the real lie group $G$ for the Hermitian symmetric pair ($G$, $K$) or merely the Lie algebra $\mathfrak{g}_{0}$. There are seven cases of Hermitian symmetric pairs which we refer to as follows:
$$\mathfrak{su}(p,q),~ \mathfrak{sp}(n, \mathbb{R}),~ \mathfrak{so}^{*}(2n),~ \mathfrak{so}(2, 2n-1),~ \mathfrak{so}(2, 2n-2),~ \mathfrak{e}_{6(-14)},~ \mathfrak{e}_{7(-25)}.$$

Let $M$ be an irreducible noncompact  Hermitian symmetric space. Let  $G$ be the identity component of the automorphism group of $M$. Let $K$ be the stabilizer of a fixed point of $M$, which is a maximal compact subgroup of $G$. Then $M$ can be identified with $G/K$ as smooth manifolds and ($G,K$) will be a Hermitian symmetric pair of noncompact type.

 We call ($G, K$) a tube type Hermitian symmetric pair when $G/K$ is a tube type Hermitian symmetric space. In our notation, the tube type Hermitian symmetric pairs are as follows:
 $$\mathfrak{su}(p,p),~ \mathfrak{sp}(n, \mathbb{R}),~ \mathfrak{so}^{*}(4n),~ \mathfrak{so}(2, 2n-1),~ \mathfrak{so}(2, 2n-2),~ \mathfrak{e}_{7(-25)}.$$

  Recall that  the irreducible unitary representations of  $G$ are in one-to-one correspondence with the irreducible unitary ($\mathfrak{g}$, $K$)-modules. And a representation of $G$ is called a highest weight
representation if its underlying ($\mathfrak{g}$, $K$)-module is a
highest weight $\mathfrak{g}$-module. Then we know a unitary highest
weight module is irreducible. When ($G,K$) is a Hermitian symmetric
pair, the unitary highest weight modules of $G$  had been classified
by Enright-Howe-Wallach \cite{EHW} and Jakobsen \cite{Jakobsen-81, Jakobsen-83}. We sometime
call such $G$-modules simply by  $\mathfrak{g}$-modules or
$\mathfrak{g}_{0}$-modules.


%

\subsection{Euclidean Jordan algebras}
The Euclidean Jordan algebras were initially introduced by P. Jordan \cite{Jordan}
for the purpose of reformulating quantum mechanics in a minimal way. Then Jordan-von Neumann-Wigner \cite{JVW} classified the simple finite dimensional Euclidean Jordan algebras: they consist of four infinity series (i.e. $\Gamma (n), \mathcal {H}_{n}(\mathbb{\mathbb{R}}),\\
 \mathcal {H}_{n}(\mathbb{\mathbb{C}}) ~and~ \mathcal {H}_{n}~(\mathbb{\mathbb{H}})$ ) and one exceptional (i.e. $\mathcal {H}_{3}~(\mathbb{\mathbb{O}})$) .

Although Euclidean Jordan algebras are abandoned by physicists
quickly, the Jordan methods have proved useful tools in a variety of
settings by mathematicians since 1950's. Some applications can be
found in McCrimmon \cite{McCrimmon-78, McCrimmon-04}.  From Koufany
\cite{Koufany} and Faraut-Kor$\mathrm{\acute{a}}$nyi \cite{Fa-Ko}, we
know the tube type hermitian symmetric spaces, irreducible tube
domains and the  conformal algebras $\mathfrak{co}(J)$ of simple
Euclidean Jordan algebras are in natural one-to-one correspondence.
So  tube type Hermitian symmetric pairs  and the conformal
algebras of simple Euclidean Jordan algebras are in one-to-one
correspondence.  In this paper, we will use the language of Euclidean
Jordan algebras to study the unitary highest weight $\mathfrak{co}(J)$-modules.

\subsection{Kepler problems and quadratic relations}
The Kepler problem is a physics problem about two bodies with an attractive force obeying the inverse square law. Mathematically it is a mechanical system with configuration space ${ \mathbb{R}}^3_{*}:={ \mathbb{R}}^3\backslash  \{ \mathbf{0}\}$ and Lagrangian
$$ L= {1\over 2}{ \mathbf{r}}'^2+{1\over r}$$
where ${\mathbf{r}}$ is a function of $t$ taking value in ${ \mathbb{R}}^3_{*}$, ${ \mathbf{r}}'$ is the velocity vector and $r$ is the length of $\mathbf{r}$. Therefore, quantum mechanically the Hamiltonian for the Kepler
problem becomes
$$
\hat{H}=-{1\over 2}\Delta-{1\over r}.
$$
where $\Delta$ is the Laplace operator on $ \mathbb{R}^{3}$.

\hspace{5cm}

The MICZ-Kepler problems are generalizations of the Kepler problems,
and they were independently discovered by McIntosh-Cisneros
\cite{MC-CI} and Zwanziger \cite{ZWANZIGER} more than thirty years ago.
The (classical) MICZ Kepler problem with magnetic charge $\mu\in {\mathbb{ R}}$ is a natural mathematical generalization of the Kepler problem,  with the Lagrangian
$$ L= {1\over 2}{\mathbf{r}}'^{2}+{1\over r} -{ \mathbf{A}}\cdot { \mathbf{r}}'-{\mu^2\over 2r^2}
$$
where ${\mathbf{A}}$ is the magnetic potential such that $ \mathbf{B}:=\nabla\times  \mathbf{A}=\mu {\mathbf{r}\over r^3}$.  Then the equation of motion is
 $${\mathbf{r}}'' = - {\mathbf{r}}'\times { \mathbf{B}}+\left({\mu^2\over r^4}-{1\over r^3}\right){ \mathbf{r}}.$$

 When Meng \cite{Meng-071, Meng-072, Meng-11, Meng-10} studied the MICZ-Kepler problem, he generalized this problem and  discovered a family of quadratic relations. Let $m$ be an integer and $m\geq 2$. We define a $(3+m)\times (3+m)$ matrix $[\eta_{\mu\nu}]=diag(1,1,-1,...,-1)$. We denote its inverse by $[\eta^{\mu\nu}]$. Let $Cl_{2,m+1}$ be
the Clifford algebra over $\mathbb{C}$ generated by $X_{\mu}'s$ subject to relations: $$X_{\mu}X_{\nu}+X_{\nu}X_{\mu}=-2\eta_{\mu\nu}.$$
We denote $M_{\mu\nu}=\frac{i}{4}(X_{\mu}X_{\nu}-X_{\nu}X_{\mu})$, then we have
\begin{equation}\label{6.4}
[M_{ab},M_{cd}]=-i(\eta_{bc}M_{ad}-\eta_{ac}M_{bd}-\eta_{bd}M_{ac}+\eta_{ad}M_{bc}). \end{equation}
 Then $\{M_{\mu\nu}\}$ generate $\mathfrak{so}(2, m+1).$

Let $(\pi, V)$ denote the unitary highest $\mathfrak{so}(2,m+1)$-module, and $\pi(\mathcal{O}):=\tilde{\mathcal{O}}$ for any
$\mathcal{O} ~in ~\mathfrak{so}(2,m+1)$. We use $[\cdot,\cdot]$  and $\{\cdot,\cdot\}$ to denote the commutator and anticommutator in a Lie algebra throughout the paper.

In Meng \cite{Meng-08}, he summarized the results in the language of representation theory:\\

 $\mathrm{\mathbf{Theorem~ \ref{Meng}}}$ {\it A (non-trivial) unitary highest weight module of $\mathfrak{so}(2,m+1)=(\mathfrak{co}(\Gamma(m)))$ has the smallest positive Gelfand-Kirillov dimension
if and only if it satisfies the following quadratic relations  on the underlying module space:
\begin{equation}\label{aa}
\{\tilde{M}_{\mu\lambda},{\tilde{M}^{\lambda}}_{\nu}\}=c\eta_{\mu\nu} ~with~\mu, \nu=-1,0,1,...,m+1,\end{equation}}
where ${\tilde{M}^{\lambda}}_{\nu}=\eta^{\lambda\kappa}\tilde{M}_{\kappa \nu}$ and $c$ is a representation-dependent real number.

\hspace{1cm}

We denote a unitary highest $\mathfrak{co}(J)$-module
by $(\pi,V)$, $\pi(\mathcal{O}):=\tilde{\mathcal{O}}$ for any
$\mathcal{O} ~in ~\mathfrak{co}(J)$.
In this paper we will show the following main theorem:\\

$\mathrm{\mathbf{Theorem~ \ref{main}}}$  {\it A (non-trivial) unitary highest weight $\mathfrak{co}(J)$-module   has the smallest positive Gelfand-Kirillov dimension
if and only if the following primary quadratic relation is
satisfied:
\begin{description}
  \item[(Q1)]
$\frac{2}{\rho}\sum\limits_{1\leq\alpha\leq
            D}\tilde{L}_{e_{\alpha}}^{2}-\tilde{L}_{e}^{2}-\frac{1}{2}\{\tilde{X}_{e},\tilde{Y}_{e}\}=-aI_{V}$ as an operator on $V$.\\
            Here D=dim$(J)$, $a=a(J, k)$(given in Remark \ref{value a}) is a nonzero constant and only depends on the highest weight $\lambda=\lambda(k)$(given in Corollary \ref{weight k}) , and $\{e_{\alpha}\}$ is an orthonormal basis for $J$. $\rho$ is the rank of $J$. $L_{e_{\alpha}},~X_{e}$ and $Y_{e}$ are generators of the conformal algebra $\mathfrak{co}(J)$.
\end{description}

}

We find that the two quadratic relations in the two theorems are equivalent when $J=\Gamma(m)$. So we call $\mathrm{(Q1)}$ a generalized quadratic relation.

\hspace{1cm}

This generalized quadratic relation is firstly introduced by Meng \cite{Meng-11(2)}.
When Meng reconstructed the various Kepler-type problems \cite{Meng-08(2), Meng-09, Meng-10(2)} in the unified language of Euclidean Jordan algebras in Ref.\cite{Meng-11(2)},  he constructed the minimal representation for the conformal algebra $\mathfrak{co}(J)$. Then he showed that this representation actually satisfied a quadratic relation, which corresponded to $k=0$ in  $\mathrm{(Q1)}$.

For a  simple complex Lie algebra $\mathfrak{g}$ not of type  $A_{n}$, Joseph \cite{Joseph-76} constructed a completely prime 2-sided ideal $J_{0}$ in the universal enveloping algebra $U(\mathfrak{g})$, whose associated variety is the closure of the minimal nilpotent orbit in $\mathfrak{g}^{*}$.  In our paper, we  find that any prime ideal $I\subseteq U(\mathfrak{g})$ equals the Joseph ideal $J_{0}$ if and only if $I$ contains  $Q_{1}'\triangleq \big(\frac{2}{\rho}\sum\limits_{1\leq\alpha\leq
            D}L_{e_{\alpha}}^{2}-L_{e}^{2}
            -\frac{1}{2}\{X_{e},Y_{e}\}+a \big)$  with $a=a(J, 0)=\frac{\rho d}{4}(1+\frac{(\rho-2)d}{4})$.

\subsection{Outline of the paper}
 In section 2, we recall basic definitions and properties about Jordan algebras. We will give their relation with tube type Hermitian symmetric pairs. In section 3, we recall the unitary highest weight modules and the classification theorem by Enright-Howe-Wallach \cite{EHW}. In section 4, we give the main theorem and some corollaries which will be used in my proof. Then in section 5, 6 and 7, we will give a case-by-case proof for our main theorem.

\section{Some basic facts in Euclidean Jordan algebras}

The facts and theorems reviewed in this section can be found in
Koecher \cite{Koecher} and Faraut-Kor$\mathrm{\acute{a}}$nyi \cite{Fa-Ko}.

A linear $\mathbf{Jordan~ algebra}$ $J$ is an(not necessarily
  associative) algebra over a field $\mathbb{F}$ whose multiplication
  satisfies the following axioms:
\begin{enumerate}
  \item $xy=yx ~(commutative~law)$.
  \item $x(x^{2}y)=x^{2}(xy)~(Jordan~identity)$.
\end{enumerate}


If we define the left multiplication in the Jordan algebra $J$ by $L_{u}(v)=uv$, then a Jordan algebra over a field $\mathbb{F}$ is just a commutative algebra over  $\mathbb{F}$ such that $[L_{u}, L_{u^2}]=0$.

Suppose $V$ is an associative algebra over a field $\mathbb{F}(Char\neq 2)$, and $A$ is a linear subspace of $V$, closed under square operations, i.e, $u\in A \Rightarrow u^{2}\in A $.
Any  such $A$  gives rise to a
Jordan algebra $A^{+}$ under a new product:
$$x\circ y:=((x+y)^2-x^2-y^2)/2.$$
For example, if  we take $V=Cl(\mathbb{R}^{n},Q)$ (i.e, the Clifford algebra of $\mathbb{R}^{n}$ equipped with a quadratic form $Q$) and $A=\mathbb{R}\oplus\mathbb{R}^{n}$, then we get a Jordan algebra $A^{+}$, denoted by $\Gamma (n)$.

 A Jordan algebra $J$ is called \textbf{special} if it can be
  realized as a Jordan subalgebra of some $A^{+}$. All other Jordan
  algebras are called \textbf{exceptional} Jordan algebras.

 A Jordan algebra $J$ is called \textbf{semi-simple} if its canonical symmetric bilinear form $\tau$(where $\tau(u,v)=the~trace~
  of~ L_{uv}$) is  non-degenerate. $J$ is called \textbf{simple} if it is
  \textbf{semi-simple} and has no nontrivial ideals.
Every semi-simple Jordan
 algebra $J$ is a direct sum of simple Jordan algebras.

A Jordan algebra $J$ over $\mathbb{R}$ is called
  \textbf{Euclidean} if its canonical symmetric bilinear form $\tau$ is positive definite.
So a \textbf{Euclidean} Jordan
  algebra is \textbf{semi-simple}.

\begin{Thm}\emph{(Jordan, von Neumann and Wigner).}
 The complete list of simple Euclidean Jordan algebras are the following:
\begin{enumerate}
  \item The algebra $\Gamma (n)=\mathbb{R}\oplus\mathbb{R}^{n} ~(n\geq2)$.
  \item The algebra $\mathcal
{H}_{n}(\mathbb{R})~(n\geq 3 ~or ~n = 1)$.
  \item The
algebra $\mathcal {H}_{n}(\mathbb{\mathbb{C}}) ~(n\geq 3)$.
  \item The algebra $\mathcal {H}_{n}~(\mathbb{\mathbb{H}})~(n\geq 3)$.
  \item The algebra $\mathcal {H}_{3}~(\mathbb{O}$).
\end{enumerate}

\end{Thm}
\begin{Rem}
$\mathcal {H}_{n}(\mathbb{F})$ denotes Hermitian $n\times n$ matrices with entries in $\mathbb{F}$.
 In 1934, Albert \cite{Albert} showed that $\mathcal {H}_{3}~(\mathbb{O}$) is not special.
\end{Rem}

\hspace{6cm}

We define the \textbf{trace} of an element in  a simple Euclidean
Jordan algebra
 $J$ as:
$tr(\lambda,\vec{ u})=2\lambda$ for $\Gamma(n)$, and
$tr(u)=trace(u)$ for the other types.

We define the \textbf{inner product} on a Euclidean Jordan algebra
$J$ by  $$<u|v> :=\frac{1}{dim(J)}\tau(u,v),$$
for all $ u, v\in J $

Every semi-simple Jordan algebra(over a field $K$ with
$\textbf{Char}=0$) has a unit element $e$.
We define the \textbf{rank} of a Euclidean Jordan algebra $J$ as
$\rho :=tr(e)$.

For a simple Euclidean Jordan algebra $J$, we have:

$$tr(u)=\frac{\rho}{dim (J)}Tr(L_{u})=\frac{\rho}{dim
 (J)}\tau(u,e).$$

Then we see $<u\mid v>=\frac{1}{\rho}tr(uv)$.

\begin{Thm}

Let $J$ be a simple Euclidean Jordan algebra of rank $\rho$ ,
$x_{0}\in J$ is non-zero and $x_{0}^{2} = tr(x_{0})x_{0}$ . Then
there is an orthogonal basis for $J\mathrm{:} \\
\{e_{11}, \cdots
, e_{\rho\rho}, e_{ij}^{\mu}\}$ with
  $\mathrm{(}$ $1 \leq i<j\leq\rho$ , $1 \leq \mu\leq d$$\mathrm{ )}$
such that

\begin{enumerate}

\item each basis vector has length $\frac{1}{\sqrt{\rho}}$;

\item $e_{ii}^{2}= e_{ii}$, $e_{ii}e_{jj}=0$ for $i\neq j$;
\item $\sum\limits_{i=1}^\rho e_{ii}=e$;

\item $(e_{jk}^{\mu})^{2}=\frac{1}{2}(e_{jj}+e_{kk})$ ,
$e_{ii}e_{ij}^{\mu}=e_{jj}e_{ij}^{\mu}=\frac{1}{2}e_{ij}^{\mu}$,
$e_{ii}e_{jk}^{\mu}=0 ~if~ i\neq j, i\neq k$ ;
\item $tre_{ii}=1$ , $tre_{ij}^{\mu}= 0$ ;

\item  $x_{0}=(trx_{0})e_{11}$ .

\end{enumerate}

\end{Thm}

\begin{Rem}The parameter $d$ in the above theorem is called the
\textbf{degree} of $J$, and $dim(J)=\rho+\frac{d\rho(\rho-1)}{2}$. The set  $\{e_{11}, \cdots
, e_{\rho\rho}\}$ is called a  $\mathbf{Jordan~frame}$.

We have the following table:

%

\begin{table}[ht]
\caption{}\label{a}
\renewcommand\arraystretch{1.5}
\noindent\[
\begin{array}{|c|c|c|c|c|c|}
\hline
J  & \Gamma(n) & \mathcal{H}_{n}(\mathbb{R})&
\mathcal{H}_{n}(\mathbb{C}) &
 \mathcal{H}_{n}(\mathbb{H})&
\mathcal{H}_{3}(\mathbb{O})      \\\hline
 \rho & 2   & n   & n  & n  & 3  \\\hline
  d     & n-1 & 1   & 2  & 4  & 8   \\\hline
\end{array}
\]
\end{table}

\end{Rem}


A $\mathbf{derivation}$ $D$ of a Jordan algebra $J$ is a linear
  transformation of $J$ such that $$D(xy)=Dx\cdot y+ x\cdot Dy,$$
  for all $ x,y\in J $. The set $\mathfrak{der}(J)$ of all derivations of $J$ is a Lie algebra with respect to the usual bracket, $[D_{1}, D_{2}]=D_{1}D_{2}-D_{2}D_{1}$.
An $\mathbf{automorphism}$ $W$ of $J$ is an invertible linear
  transformation of $J$ such that
  $$W(xy)=W(x)\cdot W(y),$$ for all $ x,y\in J $. The set $Aut(J)$ of all automorphisms of $J$ is a a Lie group since it is a  closed subgroup of $GL(J)$.


The Lie algebra of $Aut(J)$ is $\mathfrak{der}(J)$.
If $J$ is semi-simple(over $\mathbb{R}$ or $\mathbb{C}$), and $L_{x}$ is the left multiplication of $J$, then every
$D\in \mathfrak{der}(J)$ is a finite sum of $[L_{x},L_{y}]$, with $x,y \in J $.


If $J$ is a Jordan algebra, we define the $\mathbf{quadratic ~representation}$ of $J$: $$P(x)=2L_{x}^{2}-L_{x^{2}}. $$

The $\mathbf{structure~ group}$ is defined by \\
$Str(J):=\{W\in GL(J)\mid P(Wx)=WP(x)V ~\text{for all}~ x\in J ~\text{and some fixed}~ V\in GL(J)  \}.$

If we denote $S_{uv}=[L_{u},L_{v}]+L_{uv}, S_{uv}(z)=\{uvz\}$ and $S_{uv}'=S_{vu}$, then
$$[S_{uv}, S_{zw}]=S_{\{uvz\}w}-S_{z\{vuw\}}.$$
So $\mathfrak{str}(J):=span\{S_{uv} \mid u,v\in J\}$ becomes a real Lie algebra, called the \textbf{structure~algebra} of
$J$. When $J$ is a semi-simple Jordan algebra (over $\mathbb{R}$ or
$\mathbb{C}$), the Lie algebra of $Str(J)$ will equal to $\mathfrak{str}(J)$.

Let  $J$ be a simple Euclidean Jordan algebra. We denote  $\Omega$ :=\{the
  interior of $J^{2}$\}, and  $T_{\Omega}$ :=$J+i\Omega$. Let $Aut(T_{\Omega})$ be the group of holomorphic automorphisms of the tube domain $T_{\Omega}$.
 Then the Lie algebra of
$Aut(T_{\Omega})$ is \emph{$\mathfrak{co}(J)$}, called the
\textbf{conformal algebra} of $J$ with Lie bracket defined in
the next theorem.
As a vector space, $$\mathfrak{co}(J)=J\oplus \mathfrak{str}(J)\oplus J^{*}.$$
In the following we shall rewrite $u\in J$ as $X_{u}$ and $v\in J^{*}$ as $Y_{v}$.

\begin{Thm}
 (\textbf{Tits-Kantor-Koecher Construction}). Let $J$ be a simple Euclidean Jordan algebra,
then $\mathfrak{co}(J)$ := $J\oplus \mathfrak{str}(J)\oplus J^{*}$  becomes a
simple real Lie algebra with the definitions as
following:\begin{enumerate}
\item $[X_{u},X_{v}]=0$, $[Y_{u},Y_{v}]=0$;
\item $[X_{u},Y_{v}]=-2S_{uv}$;
\item $[S_{uv},X_{z}]=X_{\{uvz\}}=X_{S_{uv}(z)}$,
$[S_{uv},Y_{z}]=-Y_{\{vuz\}}=-Y_{S'_{uv} (z)}$;
\item $[S_{uv},S_{zw}]=S_{\{uvz\}w}-S_{z\{vuw\}}=S_{S_{uv}(z)w}-S_{zS'_{uv}
(w)}$.

for $u, v, z, w \in J $.
\end{enumerate}

\end{Thm}

\begin{Rem}
Actually, we have the following Table \ref{b}:

\begin{table}[ht]
\caption{}\label{b}
\renewcommand\arraystretch{1.5}
\noindent\[
\begin{array}{|c|c|c|c|}
\hline
J & \mathfrak{der}& \mathfrak{str}  &
\mathfrak{co}\\\hline
\Gamma(n)& \mathfrak{so}(n) & \mathfrak{so}(n, 1) \oplus
\mathbb{R } &\mathfrak{ so}(2, n + 1) \\\hline
\mathcal {H}_{n}(\mathbb{R}) &\mathfrak{so}(n) &
\mathfrak{sl}(n,\mathbb{R}) \oplus \mathbb{R} &
\mathfrak{sp}(n,\mathbb{R}) \\\hline
\mathcal{H}_{n}(\mathbb{C})& \mathfrak{su}(n)&
\mathfrak{sl}(n,\mathbb{C})\oplus\mathbb{ R}
&\mathfrak{su}(n,n)\\\hline
 \mathcal{H}_{n}(\mathbb{H}) &\mathfrak{sp}(n) &\mathfrak{su}^{\ast}(2n)\oplus
\mathbb{R}& \mathfrak{so}^{\ast}(4n)\\\hline
\mathcal{
H}_{3}(\mathbb{O})&  \mathfrak{f}_{4} & \mathfrak{e}_{6(-26)}
\oplus\mathbb{R} &  \mathfrak{e}_{7(-25)}\\\hline
\end{array}
\]
\end{table}

%
\end{Rem}

\section{The classification of unitary highest weight modules}

In this section we review some well-known facts and notations about the classification of unitary highest weight moduless. The details can be found in the paper Enright, Howe and Wallach
\cite{EHW}, hereafter referred to as EHW\cite{EHW}.

Let ($G, K$) be an irreducible Hermitian symmetric pair with real Lie algebras $\mathfrak{g}_{0}$ and $\mathfrak{k}_{0}$. Let $\mathfrak{g}$ and $\mathfrak{k}$ be their complexification. Then $\mathfrak{k}=\mathbb{C}H\bigoplus[\mathfrak{k}, \mathfrak{k}]$ with $ad(H)$ having eigenvalues $0, 1, -1$ on $\mathfrak{g}$.  If we denote $\mathfrak{p}^{\pm}=\{X\in \mathfrak{g}| [H, X]=\pm X\}$, then $\mathfrak{g}=\mathfrak{p}^{-}\bigoplus \mathfrak{k}\bigoplus \mathfrak{p}^{+}$.

Let $\mathfrak{h} $ be a Cartan subalgebra of both $\mathfrak{g}$ and $\mathfrak{k}$.  Let $\Delta$ denote the roots  of ($\mathfrak{g},\mathfrak{h}$) and $\Delta_{c}$ the roots  of ($\mathfrak{k},\mathfrak{h}$). Let  $\Delta_{n}$ denote the complement so that $\Delta=\Delta_{c} \cup \Delta_{n}$. We call the roots in these two sets the compact  roots and noncompact roots respectively. Let $\Delta^{+}$ denote a fixed positive root system for which $\Delta^{+}=\Delta_{c}^{+} \cup \Delta_{n}^{+}$ and $\Delta_{n}^{+}=\{\alpha\in \Delta^{+}| \mathfrak{g}_{\alpha}\subset \mathfrak{p}^{+}\}$. Let $\beta $ denote the unique maximal noncompact root of $\Delta^{+}$ and denote $\beta ^{\vee }=\frac{2\beta}{(\beta, \beta)} $. Now choose $\zeta  \in \mathfrak{h}^{*} $ so that  $\zeta $ is orthogonal to $\Delta_{c}$  and ($\zeta , \beta ^{\vee } $)=1. Let $\lambda  \in \mathfrak{h}^{*}$ be $\Delta_{c}^{+}$-dominant integral and F($\lambda  $) be the irreducible $\mathfrak{k}$-module with highest weight $\lambda $. By letting $\mathfrak{p}^{+}$ act by zero, we may consider F($\lambda  $) as a module of $\mathfrak{q}=\mathfrak{k}\oplus \mathfrak{p}^{+}$. Then we define:
$$N(\lambda )=U(\mathfrak{g})\otimes_{U(\mathfrak{q})}F(\lambda ).$$

Let $L(\lambda)$ denote the irreducible quotient of $ N(\lambda )$. If $L(\lambda)$  is unitarizable(i.e it is equivalent to the $\mathfrak{g}$-module of $K$-finite vectors in a unitary representation of $G$, or it is a unitary ($\mathfrak{g}$, $K$)-module), then $\lambda=\lambda_{0}+z \zeta $, with $\lambda_{0} \in \mathfrak{h}^{*} $ such that ($\lambda_{0}  + \rho, \beta  $)=0, and $z\in \mathbb{R}$. For a fixed $\lambda_{0}$, the set of all $z\in \mathbb{R}$ with $L(\lambda_{0}+z \zeta )$ a unitary $\mathfrak{g}$-module takes the  following form:
\begin{center}
\setlength{\unitlength}{1mm}
\begin{picture}(180,4)

\line(1,0){60}
\put(0,-1){$\bullet $}
\put(-1,-5){$A(\lambda_{0})$}
\put(10,-1){$\bullet $}
\put(20,-1){$\bullet $}

\put(27,-0.5){$...~~...$ }
\put(45,-1){$\bullet $}
\put(55,-1){$\bullet $}
\put(54,-5){$B(\lambda_{0})$}

\end{picture}

\end{center}

\hspace{1cm}

Let $Z(\lambda_{0})=\{z\in \mathbb{R}\mid L(\lambda_{0}+z \zeta ) ~is  ~unitarizable~ \}$. Let $Z(\lambda_{0})_{r}$ denote the subset of $Z(\lambda_{0})$ for which $ N(\lambda )$ is reducible. We call
$Z(\lambda_{0})_{r}$ the unitary reduction points on the line. These points correspond to the value of z on the above line which are equally spaced from $A(\lambda_{0})$ to
$B(\lambda_{0})$. $A(\lambda_{0})$ is called the first reduction point, and $B(\lambda_{0})$ is called the last reduction point. The $\mathbf{reduction~ level}$ $r(\lambda_{0})$ of  $\lambda_{0}$ is the number of points in $Z(\lambda_{0})_{r}$. We list these reduction points by $$A(\lambda_{0})=z_{r(\lambda_{0})}^{\lambda_{0}}<z_{r(\lambda_{0})-1}^{\lambda_{0}}<...<z_{1}^{\lambda_{0}}=B(\lambda_{0}).$$ From now on we simply call $L(\lambda)$ a unitary highest weight $\mathfrak{g}_{0}$-module since it is a module of the hermitian symmetric pair $\mathfrak{g}_{0}$.

We define $C(\lambda_{0})=(B(\lambda_{0})-A(\lambda_{0}))/(r(\lambda_{0})-1)$. Let $r$ be the split rank of  $\mathfrak{g}_{0}$. Then $r$ equals the reduction level of the weight $\lambda_{0}=-(\rho,\beta^{\vee})\zeta$. And we have the following table from EHW\cite{EHW}:

\begin{table}[ht]
\caption{}\label{d}
\resizebox{14cm}{!}{$
\begin{array}
{|c|c|c|c|c|c|c|c|}\hline
 \mathfrak{g}_{0}& \mathfrak{su}(p,q) & \mathfrak{sp}(n,\mathbb{R}) &\mathfrak{so}^{*}(2n)& \mathfrak{so}(2,2n-2) & \mathfrak{so}(2,2n-1) & \mathfrak{e}_{6(-14)}& \mathfrak{e}_{7(-25)}   \\\hline
  C(\lambda_{0})               & 1 & 1/2 & 2& n-2 & n-3/2 & 3 & 4 \\\hline
r& min\{p,q\} &n&  [n/2]  &2&2&2&3 \\\hline
(\rho,\beta^{\vee})& p+q-1&n &2n-3&2n-3&2n-2&11 &17 \\\hline
\end{array}
$}
\end{table}

So we can denote $C(\lambda_{0})$ by $C$ since it is independent of $\lambda_{0}$.

\hspace{4cm}

We have the following property:
\begin{Lem}
Suppose $J$ is one of the 5 simple Euclidean Jordan algebras, then
\begin{enumerate}
  \item $r=\rho$, i.e., the split rank of  $co(J)$ equals the rank of $J$;
  \item $2C=d$, \text{here $d$ is the degree of $J$}.
\end{enumerate}

\end{Lem}
%

%
%
Let $\mathfrak{g}$ be a complex Lie algebra and $U(\mathfrak{g})$ be the enveloping
algebra of $\mathfrak{g}$. We denote by $U_{n}(\mathfrak{g})$ the finite dimensional subspace of $U(\mathfrak{g})$, which are spanned by products of at most n elements in $\mathfrak{g}$. Then $\{U_{n}(\mathfrak{g})\}_{n=0}^{\infty }$ is the natural filtration of $U(\mathfrak{g})$. By the Poincar$\mathrm{\acute{e}}$-Birkhoff-Witt theorem, we can identify the associated graded ring $\mathrm{gr} U(\mathfrak{g})=\oplus_{n=0}^{\infty }U_{n}(\mathfrak{g})/U_{n-1}(\mathfrak{g})$
 with the symmetric algebra $S(\mathfrak{g}) $.

Let M be a finitely generated $U(\mathfrak{g})$-module. We can take a finite-dimensional subspace
$M_{0}$ of M such that $M=U(\mathfrak{g})M_{0}$. Define $M_{n}=U_{n}(\mathfrak{g})M_{0}(n=1,2,... )$.
Then $\{M_{n}\}$ gives an increasing filtration of M. We get a
finitely generated graded $S(\mathfrak{g})$-module $\mathrm{gr}M=\oplus_{n=0}^{\infty }\mathrm{gr}_{n}M, ~\mathrm{gr}_{n}M=M_{n}/M_{n-1}.$

Let $\mathrm{Ann}(L(\lambda))\subseteq U(\mathfrak{g}) $ be the annihilator of the unitary highest weight module $L(\lambda)$.
Then its graded ideal is: $$\mathrm{gr}(\mathrm{Ann}(L(\lambda)))=\mathrm{Ann}_{S(\mathfrak{g})}(\mathrm{gr}L(\lambda))=\{D\in S(\mathfrak{g})\mid Dw=0 ~for~all~w\in \mathrm{gr}(L(\lambda))\}.$$    The $\mathbf{associated ~variety}$ of the ($\mathfrak{g}$, $K$)-module $L(\lambda)$ is defined to be:
$$\mathcal {V}(L(\lambda))=\mathcal {V}(\mathrm{gr}(\mathrm{Ann}(L(\lambda))))=\{X\in \mathfrak{g}^{*}\mid D(X)=0 ~for~all~D\in \mathrm{gr}(\mathrm{Ann}(L(\lambda)))\}.$$
Here $S(\mathfrak{g}) $ is viewed as
the polynomial ring over $\mathfrak{g}^{*}$ through the Killing form of $\mathfrak{g}$ . In the language of algebraic geometry, $\mathcal {V}(L(\lambda))$ is just the set of all prime ideals containing $Ann_{S(\mathfrak{g})}\mathrm{gr}(L(\lambda))$, a closed subvariety of the affine algebraic variety Spec$S(\mathfrak{g})$. It is also a $K$-invariant subvariety of $\mathfrak{g}^{*}$(or of $\mathfrak{(g/k)}^{*}$) since the compatibility of the $K$ action and the module structure on $L(\lambda)$. In fact  $\mathcal {V}(L(\lambda))$ is equal to the closure of some nilpotent $K_{\mathbb{C}}$-orbit, also equal to the union of a finite number of nilpotent $K_{\mathbb{C}}$-orbits  by Vogan \cite{Vogan-91}.

We define the $\mathbf{Gelfand}$-$\mathbf{Kirillov ~dimension}$ of $L(\lambda)$ by:

$$GKdim(L(\lambda))=dim\mathcal {V}(L(\lambda)). $$

For a  simple complex Lie algebra $\mathfrak{g}$ not of type  $A_{n}$, Joseph \cite{Joseph-76} constructed a unique completely prime 2-sided ideal $J_{0}$ in the universal enveloping algebra $U(\mathfrak{g})$, with the property that its associated variety $\mathcal{V}(\mathrm{gr}(J_{0}))$ is equal to the closure of the minimal nilpotent orbit in $\mathfrak{g}^{*}$. This ideal is callled the $\mathbf{Joseph~Ideal}.$ However, it was noticed by Savin \cite{Savin} that there is a gap in the proof of Lemma 8.8 in Joseph\cite{Joseph-76}. Gan and Savin\cite{Ga-Sa} proved the uniqueness property of the Joseph ideal.

A unitary highest weight $\mathfrak{g}_{0}$-module $L(\lambda)$ is called $\mathbf{minimal}$ if its annihilator is the Joseph ideal. So it must have the smallest positive Gelfand-Kirillov dimension.
From refs \cite{EN-HU, En-Wi, Joseph-92}, we have the following lemma.

\begin{Lem}
 A unitary highest weight $\mathfrak{g}_{0}$-module $L(\lambda)=L(\lambda_{0}+z\zeta)$ has the smallest positive Gelfand-Kirillov dimension  if and only if $z=(\rho,\beta^{\vee})-C$, and $L(\lambda_{0}+z\zeta)$ has the zero Gelfand-Kirillov dimension  if and only if $z=(\rho,\beta^{\vee})$.

\end{Lem}

Following the notations in Bourbaki \cite{Bourbaki}, we can write the highest weight $\lambda$ as a tuple of real numbers. Then we have the following corollary.

\begin{Cor}\label{weight k}
The unitary highest weight $\mathfrak{co}(J)$-module $L(\lambda)=L(\lambda_{0}+z\zeta)$ has the smallest positive Gelfand-Kirillov dimension  if and only if: \\
$\lambda=\left\{
  \begin{array}{ll}
  (-(n+k-\frac{1}{2}),k,...,k), k=0 ~\text{\emph{or}} ~1/2,  & \text{\emph{if}} ~J=\Gamma(2n)\\
  (-(n+|k|-1),|k|,...,|k|,k), k ~\text{\emph{is a half integer}},  & \text{\emph{if}} ~J=\Gamma(2n-1)\\

  -(\frac{1}{2},\frac{1}{2},...,\frac{1}{2},\frac{1}{2}+k), k=0,1. & \text{\emph{if}}  ~ J=\mathcal{H}_{n}(\mathbb{R})\\
(\underbrace{-\frac{n+k}{2n},\dots,-\frac{n+k}{2n}}_{n},\frac{n-k}{2n}+k,\frac{n-k}{2n},...,\frac{n-k}{2n}),k=0,1,... \\
or~ (-\frac{n-k}{2n},...,-\frac{n-k}{2n},-\frac{n-k}{2n}-k,\underbrace{\frac{n+k}{2n},...,\frac{n+k}{2n}}_{n}), k=0,1,...& \text{\emph{if}} ~J=\mathcal{H}_{n}(\mathbb{C}) \\
-(1,...,1,k+1),k=0,1,... &\text{\emph{if}}  ~ J=\mathcal{H}_{n}(\mathbb{H})\\
(0,0,0,0,0,-4,2,-2),      & \text{\emph{if}}  ~J=\mathcal{H}_{3}(\mathbb{O})\\
  \end{array}
\right.$

~\\

$L(\lambda_{0}+z\zeta)$ has the zero GK dimension  if and only if $\lambda$ equals zero weight.

%

\end{Cor}

{\bf Proof.}
 From the above lemma,  a unitary highest weight $\mathfrak{g}_{0}$-module $L(\lambda)=L(\lambda_{0}+z\zeta)$ has the smallest positive Gelfand-Kirillov dimension means $z=(\rho,\beta^{\vee})-C$. Then case by case from  EHW\cite{EHW},  we can compute all possible $\lambda_{0}$  for this $z$. $\square$

\hspace{1cm}

From Hilgert-Kobayashi-M$\mathrm{\mathrm{\ddot{o}}}$llers \cite{HKM} and Torasso \cite{Torasso}, we have the following lemma about minimal representations.

\begin{Lem}
A unitary highest weight $\mathfrak{co}(J)$-module $L(\lambda)=L(\lambda_{0}+z\zeta)$ is minimal  if and only if: \\
$\lambda=\left\{
  \begin{array}{ll}
  (-(n-\frac{1}{2}),0,...,0),   & \text{\emph{if}} ~J=\Gamma(2n)\\
  (-(n-1),0,...,0),   & \text{\emph{if}} ~J=\Gamma(2n-1)\\
  -(\frac{1}{2},\frac{1}{2},...,\frac{1}{2},\frac{1}{2}+k), k=0,1. & \text{\emph{if}}  ~ J=\mathcal{H}_{n}(\mathbb{R})\\
(-1,...,-1),           &\text{\emph{if}}  ~ J=\mathcal{H}_{n}(\mathbb{H})\\
(0,0,0,0,0,-4,2,-2),      & \text{\emph{if}}  ~J=\mathcal{H}_{3}(\mathbb{O})\\
  \end{array}
\right.$

~\\
In other words, $\lambda=-C\zeta$ except for the odd part of the metaplectic representation.

\end{Lem}
The unitary highest $\mathfrak{g}_{0}$-module $L(\lambda)=L(-C\zeta)$ is called the first Wallach representation.

\section{Main theorem}
We denote a unitary highest $\mathfrak{co}(J)$-module $L(\lambda)=L(\lambda_{0}+z\zeta)$
by $(\pi,V)$, $\pi(\mathcal{O}):=\tilde{\mathcal{O}}$ for any
$\mathcal{O} ~in ~\mathfrak{co}(J)$.
\begin{Thm}\label{main}
A (non-trivial) unitary highest weight $\mathfrak{co}(J)$-module $L(\lambda)=L(\lambda_{0}+z\zeta)$  has the smallest positive Gelfand-Kirillov dimension
if and only if the following primary quadratic relation is
satisfied:
\begin{description}
  \item[(Q1)]
$\frac{2}{\rho}\sum\limits_{1\leq\alpha\leq
            D}\tilde{L}_{e_{\alpha}}^{2}-\tilde{L}_{e}^{2}-\frac{1}{2}\{\tilde{X}_{e},\tilde{Y}_{e}\}=-aI_{V}$ as an operator on $V$.\\
            Here D=dim$(J)$, $a=a(J, k)$(the explicit values will be given in Remark \ref{value a}) is a nonzero constant and only depends on the highest weight $\lambda=\lambda(k)$ given in Corollary \ref{weight k}, and $\{e_{\alpha}\}$ is an orthonormal basis for $J$.
\end{description}

\end{Thm}

\begin{Cor}\label{Q2-4}
From the quadratic relation $\mathrm{(Q1)}$, we can get the following secondary quadratic relations:
\begin{center}
\begin{description}
          \item[(Q2)] $\sum\limits_{1\leq\alpha\leq D}\tilde{X}_{e_{\alpha}}^{2}=\rho
            \tilde{X}_{e}^{2}$.
            \item[(Q3)]  $\sum\limits_{1\leq\alpha\leq D}\{\tilde{X}_{e_{\alpha}},\tilde{X}_{e_{\alpha}u}\}=\rho\{\tilde{X}_{e},\tilde{X}_{u}\}$,
for any $u\in J.$

           \item[(Q4)]  $\frac{2}{\rho}\sum\limits_{1\leq\alpha\leq D}\{\tilde{L}_{e_{\alpha}},\tilde{L}_{e_{\alpha}u^{2}-u(ue_{\alpha})}\}+\frac{4}{\rho}\sum\limits_{1\leq\alpha \leq
            D}[\tilde{L}_{u},\tilde{L}_{e_{\alpha}}]^{2}
            +\{\tilde{X}_{u},\tilde{Y}_{u}\}
            -\frac{1}{2}\{\tilde{X}_{u^{2}},\tilde{Y}_{e}\}
            -\frac{1}{2}\{\tilde{X}_{e},\tilde{Y}_{u^{2}}\}=0$,
~for any $u\in J.$
 \end{description}
\end{center}
\end{Cor}
{\bf Proof of the corollary.}
Meng \cite{Meng-11(2)} has proved that the quadratic relation (Q1) implies (Q2)
for all $J$.\\

We compute $[(Q1), \tilde{X}_{e}]$ and get:
\begin{equation}\label{5.1}
\sum\limits_{1\leq \alpha \leq D}\{ \tilde{X}_{e_{\alpha}},\tilde{L}_{e_{\alpha}}\}=\rho\{ \tilde{X}_{e},\tilde{L}_{e}\}.\end{equation}

Then $[(\ref{5.1}), \tilde{X}_{u}]$ implies (Q3).

We compute $[(Q1), \tilde{L}_{u}]$ and get:

\begin{equation}\label{5.2}
\frac{2}{\rho}\sum\limits_{1\leq\alpha \leq        D}\{\tilde{L}_{e_{\alpha}},[\tilde{L}_{e_{\alpha}},\tilde{L}_{u}]\}
+\frac{1}{2}(\{\tilde{X}_{u}, \tilde{Y}_{e} \}-\{\tilde{X}_{e}, \tilde{Y}_{u} \})=0.
\end{equation}

Then $[(\ref{5.2}), \tilde{L}_{u}]$ implies (Q4). $\Box$

\begin{Cor}\label{minimal}
Suppose $\mathfrak{g}_{0}=\mathfrak{co}(J)\neq \mathfrak{su}(n,n)$ (i.e. not of type A) and $(\pi, V)$ is a (non-trivial) unitary highest weight $\mathfrak{co}(J)$-module. Then $(\pi, V)$  is minimal if and only if the following quadratic relation is
satisfied:
\begin{equation}
 \frac{2}{\rho}\sum\limits_{1\leq\alpha\leq
            D}\tilde{L}_{e_{\alpha}}^{2}-\tilde{L}_{e}^{2}-\frac{1}{2}\{\tilde{X}_{e},\tilde{Y}_{e}\}=-aI_{V} \text{~as~an~operator~on~} V.
           \end{equation}
 Here $D=dim(J)$, $a=\frac{\rho d}{4}(1+\frac{(\rho-2)d}{4})$ is a nonzero constant and equal to $a(J, 0)$ in our main theorem, and $\{e_{\alpha}\}$ is an orthonormal basis for $J$.

\end{Cor}
\hspace{4cm}

\begin{Rem}The quadratic relation in the above corollary is a special case of $(\mathrm{Q1})$. Actually all highest weights in a minimal $\mathfrak{co}(J)$-module correspond to $a=a(J,0)$ in Corollary \ref{value a}. From Meng\cite{Meng-11(2)}, we know $a(J,0)=\frac{\rho d}{4}(1+\frac{(\rho-2)d}{4})$.

For a  simple complex Lie algebra $\mathfrak{g}=\mathfrak{co}(J)_{\mathbb{C}}$ not of type  $A_{n}$, we know that the Joseph ideal  $J_{0}$ in $U(\mathfrak{g})$ is the unique completely prime  ideal  whose associated variety is the closure of the minimal nilpotent orbit in $\mathfrak{g}^{*}$.  From the definition of minimality and our Corollary \ref{minimal}, we must have: \\

{\it the annihilator ideal $Ann(L(\lambda))$ in the universal enveloping algebra $U(\mathfrak{co}(J))_{\mathbb{C}})$ is equal to the Joseph ideal $J_{0}$ if and only if $Ann(L(\lambda))$ contains $$Q_{1}'\triangleq \big(\frac{2}{\rho}\sum\limits_{1\leq\alpha\leq
            D}L_{e_{\alpha}}^{2}-L_{e}^{2}
            -\frac{1}{2}\{X_{e},Y_{e}\}+\frac{\rho d}{4}(1+\frac{(\rho-2)d}{4}) \big).$$}

~\\
Then any prime ideal $I\subseteq U(\mathfrak{g})$ equals the Joseph ideal $J_{0}$ if and only if $I$ contains  $Q_{1}'\triangleq \big(\frac{2}{\rho}\sum\limits_{1\leq\alpha\leq
            D}L_{e_{\alpha}}^{2}-L_{e}^{2}
            -\frac{1}{2}\{X_{e},Y_{e}\}+a \big)$  with $a=a(J, 0)=\frac{\rho d}{4}(1+\frac{(\rho-2)d}{4})$.

\hspace{1cm}

In ref.\cite{Meng-11(2)},  Meng constructed the models of some generalized Kepler problems. He found that:
a minimal unitary highest weight $\mathfrak{co}(J)$-module can be realized by a $L^{2}$-model, i.e. the Hilbert space of bound states of a model from a generalized Kepler problem he defined in his paper.  Then Meng showed: the quadratic relation $(Q1)$ with $a=a(J, 0)$ is in fact in the annihilator of the minimal $\mathfrak{co}(J)$-module he defined.

\end{Rem}

\hspace{1cm}

Before we  give a proof for our main theorem, we define some notations which will be used throughout the rest of this paper.
We take the same notation with Meng\cite{Meng-11(2)}. Denote
$$E_{u}^{\pm}=iL_{u}\mp\frac{1}{2}(X_{u}-Y_{u}),~~ h_{u}=-i(X_{u}+Y_{u})$$
for any $u\in J$.

Then we have the following property according to the TKK commutation relations:

\begin{Lem}
~\\
\begin{enumerate}
  \item $[h_{u}, E_{v}^{\pm}]=\pm2E_{uv}^{\pm}$;
  \item $[E_{u}^{+}, E_{v}^{-}]=-h_{uv}-2[L_{u}, L_{v}]$;
  \item $[E_{u}^{+}, E_{v}^{+}]=[E_{u}^{-}, E_{v}^{-}]=0$;
   \item $[h_{u}, h_{v}]=4[L_{u}, L_{v}]$.
\end{enumerate}

\end{Lem}

\hspace{1cm}

\section{Proof of the main theorem-the case when $J=\Gamma(m)(m\geq2)$ }

Actually this case had been proved by
Meng \cite{Meng-08}.

Let $m$ be an integer and $m\geq 2$. We define a $(3+m)\times (3+m)$ matrix $[\eta_{\mu\nu}]=diag(1,1,-1,...,-1)$. We denote its inverse by $[\eta^{\mu\nu}]$. Let $Cl_{2,m+1}$ be
the Clifford algebra over $\mathbb{C}$ generated by $X_{\mu}'s$ subject to relations: $$X_{\mu}X_{\nu}+X_{\nu}X_{\mu}=-2\eta_{\mu\nu}.$$
We denote $M_{\mu\nu}=\frac{i}{4}(X_{\mu}X_{\nu}-X_{\nu}X_{\mu})$, then we have
\begin{equation}\label{6.4}
[M_{ab},M_{cd}]=-i(\eta_{bc}M_{ad}-\eta_{ac}M_{bd}-\eta_{bd}M_{ac}+\eta_{ad}M_{bc}). \end{equation}
 Then $\{M_{\mu\nu}\}$ generate $\mathfrak{so}(2, m+1).$

Let $(\pi, V)$ denote the unitary highest $\mathfrak{so}(2,m+1)$-module, and $\pi(\mathcal{O}):=\tilde{\mathcal{O}}$ for any
$\mathcal{O} ~in ~\mathfrak{so}(2,m+1)$.
In \cite{Meng-08},  Meng proved the following theorem.
\begin{Thm}$(\mathrm{Meng}).$\label{Meng}
A  (non-trivial) unitary highest weight module of $\mathfrak{so}(2,m+1)=(\mathfrak{co}(\Gamma(m)))$ has the smallest positive Gelfand-Kirillov dimension
if and only if it satisfies the following quadratic relations on the underlying module space:
\begin{equation}\label{aa}
\{\tilde{M}_{\mu\lambda},{\tilde{M}^{\lambda}}_{\nu}\}=c\eta_{\mu\nu} ~with~\mu, \nu=-1,0,1,...,m+1,\end{equation}
where ${\tilde{M}^{\lambda}}_{\nu}=\eta^{\lambda\kappa}\tilde{M}_{\kappa \nu}$ and $c$ is a representation-dependent real number.
\end{Thm}
Then from (\ref{6.4})  we find that:
\begin{Lem}

The quadratic relation in Meng's theorem is equivalent to the following relation:
\begin{equation}
\{M_{0\lambda},{M^{\lambda}}_{0}\}=c\eta_{00}=c,\end{equation} plus the  commutation relations between the M's in (\ref{6.4}).
\end{Lem}

\

 The Jordan algebra $\Gamma(m)$ has an orthogonal basis:\\
  $\{e_{11}=(\frac{1}{2},\frac{1}{2},0,...,0),
e_{22}=(\frac{1}{2},-\frac{1}{2},0,...,0), e_{12}^{1}=(0,0,\frac{1}{\sqrt{2}},0,...,0)=\frac{1}{\sqrt{2}}e_{3},\\ e_{12}^{\mu}=(0,0,...,0,\frac{1}{\sqrt{2}},0,...,0)=\frac{1}{\sqrt{2}}e_{\mu+2}| 1\leq\mu \leq d=m-1\}$.

\hspace{5cm}

We define the following linear transformation:
\begin{equation}
\begin{split}
\varphi: ~ co(\Gamma(m))&\rightarrow so(2,m+1)\\
             X_{e}&\mapsto -i(M_{-1,0}+M_{0,m+1}),\\
             Y_{e}&\mapsto -i(M_{-1, 0}-M_{0,m+1}),\\
            L_{(\lambda,\vec{u}_{0})}&\mapsto -i(-\lambda
                        M_{-1,m+1}+\sum\limits_{1\leq i\leq m}u_{i}M_{0,i}), \\
                        &\text{here}~\vec{u}_{0}=(u_{1},u_{2},...,u_{m}).
\end{split}
\end{equation}

By using the TKK commutation relations, we can get
$$\varphi(X_{e_{11}})=-\frac{i}{2}(M_{-1,0}+M_{-1,1}+M_{0,m+1}+M_{1,m+1}),$$
$$\varphi(X_{e_{22}})=-\frac{i}{2}(M_{-1,0}-M_{-1,1}+M_{0,m+1}-M_{1,m+1}),$$
$$\varphi(Y_{e_{11}})=-\frac{i}{2}(M_{-1,0}-M_{-1,1}-M_{0,m+1}+M_{1,m+1}),$$
$$\varphi(Y_{e_{22}})=-\frac{i}{2}(M_{-1,0}+M_{-1,1}-M_{0,m+1}-M_{1,m+1}),$$
$$\varphi(X_{e_{12}^{\mu}})=-\frac{i}{\sqrt{2}}(M_{-1,\mu+1}+M_{\mu+1,m+1}),$$
$$\varphi(Y_{e_{12}^{\mu}})=-\frac{i}{\sqrt{2}}(-M_{-1,\mu+1}+M_{\mu+1,m+1}).$$
Then we can easily check that $\varphi$ is an isomorphism by computing its inverse $\varphi^{-1}$ and $\varphi^{-1}$ satisfies the commutators.

\hspace{6cm}

 From this linear transformation we find that:
 $$\{M_{0\lambda},{M^{\lambda}}_{0}\}=c\eta_{00}=c:=-a$$ is equivalent to

 $$\frac{2}{\rho}\sum\limits_{1\leq\alpha\leq
            D}L_{e_{\alpha}}^{2}-L_{e}^{2}
            -\frac{1}{2}\{X_{e},Y_{e}\}=-a.$$
And the commutation relations in (\ref{6.4}) is equivalent to the TKK commutation relations.

\hspace{6cm}

The TKK commutation relations is naturally satisfied by such unitary highest weight $\mathfrak{co}(J)$-modules from the definition of Lie algebra modules.  So Meng's theorem is equivalent to:

\begin{Thm}
A  (non-trivial) unitary highest weight module of $\mathfrak{so}(2,m+1)=(\mathfrak{co}(\Gamma(m)))$ has the smallest positive Gelfand-Kirillov dimension
if and only if the following quadratic relation is
satisfied as operators on the underlining module space $V$:
$$\frac{2}{\rho}\sum\limits_{1\leq\alpha\leq
            D}\tilde{L}_{e_{\alpha}}^{2}-\tilde{L}_{e}^{2}
            -\frac{1}{2}\{\tilde{X}_{e},\tilde{Y}_{e}\}=-a.$$
\end{Thm}

%
\vspace{1cm}

\section{Proof of the main theorem-the hermitian cases (Part I)  }
\hspace{1cm}

For the remaining hermitian cases, firstly we suppose a unitary highest weight $\mathfrak{co}(J)$-module  satisfies a quadratic relation $(\mathrm{Q}1)$ in the main theorem. We will show that it will have the smallest positive Gelfand-Kirillov dimension by  computing its highest weight in a case-by-case way. Then in the next section, we will show the other direction of the main theorem.

The idea is very simple. Actually, if we denote $Q_{1}\triangleq \big(\frac{2}{\rho}\sum\limits_{1\leq\alpha\leq
            D}\tilde{L}_{e_{\alpha}}^{2}-\tilde{L}_{e}^{2}
            -\frac{1}{2}\{\tilde{X}_{e},\tilde{Y}_{e}\}+a \big)$, then the quadratic relation $(\mathbf{Q1})$ is equivalent to  $(Q_{1})v_{\lambda}=0$ and $[(Q_{1}), \pi (E^{n_{1}}_{\beta_{1}}\cdot\cdot\cdot E^{n_{k}}_{\beta_{k}})]v_{\lambda}=0$, where $\beta_{i}\in \Delta$ is a root and $v_{\lambda}$ is any highest weight vector in $V$. From theses equations, we can work out the highest weight $\lambda$. Mainly we will use the equations $(\mathbf{Q2})$, $(\mathbf{Q3})$ and $(\mathbf{Q4})$ to compute the highest weight $\lambda$.

\subsection{The case when $J=\mathcal {H}_{n}(\mathbb{R})(n\geq3)$}

\hspace{0.5cm}

We will follow the approach in Meng \cite{Meng-08}. The idea
is to construct a convenient Cartan basis and then rewrite the
quadratic relations in terms of the Cartan basis elements.

In this case, $\mathfrak{co}(J)=\mathfrak{sp}(n,\mathbb{R})$. The
positive roots are $e_{i}-e_{j}$(compact roots),
$e_{i}+e_{j}$(noncompact roots) with $i<j$ and $2e_{i}$(noncompact
roots) with $1\leq i\leq n$. We
choose the following Cartan basis for $\mathfrak{sp}(n,\mathbb{R})$:

\begin{align*}
&H_{2e_{i}}=h_{e_{ii}}=-i(X_{e_{ii}}+Y_{e_{ii}})&  H_{e_{i}\pm e_{j}}&=h_{e_{ii}}\pm h_{e_{jj}}\\
&E_{2e_{i}}=E_{e_{ii}}^{+}=iL_{e_{ii}}-\frac{1}{2}(X_{e_{ii}}-Y_{e_{ii}})&  E_{-2e_{i}}&=E_{e_{ii}}^{-}=iL_{e_{ii}}+\frac{1}{2}(X_{e_{ii}}-Y_{e_{ii}})\\
&E_{e_{i}+e_{j}}=\sqrt{2}E_{e_{ij}}^{+} & E_{-e_{i}-e_{j}}&=\sqrt{2}E_{e_{ij}}^{-} \\
&E_{\pm(e_{i}-e_{j})}=\frac{1}{\sqrt{2}}(h_{e_{ij}}\pm 4[L_{e_{ij}},L_{e_{jj}}])
\end{align*}

Then we get:
 \begin{align*}
      L_{e_{ii}}&=-\frac{i}{2}(E_{2e_{i}}+E_{-2e_{i}}) & L_{e_{ij}}&=-\frac{i}{2\sqrt{2}}(E_{e_{i}+e_{j}}+E_{-e_{i}-e_{j}})\\
      X_{e_{ii}}&=\frac{1}{2}(E_{-2e_{i}}-E_{2e_{i}}+iH_{2e_{i}})& Y_{e_{ii}}&=-\frac{1}{2}(E_{-2e_{i}}-E_{2e_{i}}-iH_{2e_{i}})\\
      X_{e_{ij}}&=\frac{1}{2\sqrt{2}}(iE_{e_{i}-e_{j}}+iE_{-e_{i}+e_{j}}+E_{-e_{i}-e_{j}}-E_{e_{i}+e_{j}})&&\\
      Y_{e_{ij}}&=\frac{1}{2\sqrt{2}}(iE_{e_{i}-e_{j}}+iE_{-e_{i}+e_{j}}-E_{-e_{i}-e_{j}}+E_{e_{i}+e_{j}})&&\\
  \end{align*}

Let
($\pi$, V) be the unitary  highest weight module of $\mathfrak{co}(J)=\mathfrak{sp}(n,\mathbb{R})$ with highest weight $\lambda=(\lambda_{1},\lambda_{2},...,\lambda_{n})$.
Then for any highest weight vector $v \in V$, we have
$\pi (H_{\alpha})v=\tilde{H}_{\alpha}v=\lambda(H_{\alpha})v=\frac{2(\lambda,\alpha)}{(\alpha,\alpha)}v$
and $\pi (E_{\alpha})v=0$ if $\alpha$ is a positive root.

We let  the quadratic relations (Q1) and (Q2) act on any highest weight vector, then we get:

\begin{equation}\label{4.1}
(\sum\limits_{1\leq i\leq n}\lambda_{i})^{2}+(n+1)(\sum\limits_{1\leq i\leq n}\lambda_{i})+4a=0.
\end{equation}

and

 \begin{equation}\label{4.2}
  (\sum\limits_{1\leq i\leq n}\lambda_{i})^{2}-\sum\limits_{1\leq i\leq
     n}\lambda_{i}^{2}\\
     +\sum\limits_{1\leq i<j\leq n}\lambda_{j}=0.
\end{equation}

By the unitarity, we can get:
 \begin{equation}\label{6.8}
\lambda_{n}\leq
\lambda_{n-1}\leq ...\leq \lambda_{1}\leq 0. \end{equation}

%
%
%
%
%
%
%
%
%
%
%
%
%

From the quadratic relation (Q3), we can get
 \begin{eqnarray}\label{e11}
    \left\{\begin{array}{l}
 \sum\limits_{1< j\leq n}\lambda_{j}+2\lambda_{1}(\sum\limits_{1< j\leq n}\lambda_{j})=0, \text{ here we take}~  u=\sqrt{\rho}e_{11}.\\
~\\
\sum\limits_{ i< j\leq n}\lambda_{j}+\sum\limits_{ 1\leq k< i}\lambda_{i}+2\lambda_{i}(\sum\limits_{p\neq i}\lambda_{p})=0, \\
\text{ here we take}~  u=\sqrt{\rho}e_{ii}, 2\leq i \leq n-1.\\
~\\
\sum\limits_{1\leq k< n}\lambda_{n}+2\lambda_{n}(\sum\limits_{1\leq j< n}\lambda_{j})=0, \text{ here we take}~  u=\sqrt{\rho}e_{nn}.\\

\end{array}
\right.
\end{eqnarray}

Since $\lambda$ is not a zero weight, then from (\ref{e11}) and (\ref{6.8}),
 we must have
$\lambda_{1}=\lambda_{2}=...=\lambda_{n-1}=-\frac{1}{2}\geq\lambda_{n}
$.

From the quadratic relation (Q4), we can get
\begin{equation}\label{q4e1n}
(2\lambda_{1}+2\lambda_{n}+1)\sum\limits_{1\leq i\leq n}\lambda_{i}+(3+n)\lambda_{n}-3\lambda_{1}=0, \text{ here we take}~  u=\sqrt{\rho}e_{1n}.
\end{equation}

Then we can get
$\lambda_{n}=-\frac{1}{2}-k$, with $k=0,1$.





~\\
So
$\lambda=\left\{
          \begin{array}{ll}
            -(\frac{1}{2},\frac{1}{2},...,\frac{1}{2}), & \hbox{if k=0} \\
            -(\frac{1}{2},\frac{1}{2},...,\frac{1}{2},\frac{1}{2}+1), & \hbox{if k=1}
          \end{array}
        \right.$

\hspace{2cm}



\subsection{The case when $J=\mathcal {H}_{n}(\mathbb{C})(n\geq3)$}
\hspace{1cm}

In this case, $\mathfrak{co}(J)=\mathfrak{su}(n,n)$. The compact
positive roots are $e_{i}-e_{j}$ with $1\leq i<j\leq n$ or $n+1\leq
i<j\leq 2n$, and the noncompact roots are $e_{i}-e_{j}$ with $1\leq
i\leq n$ and $2n\geq j\geq n+1$. Let
$\{\sqrt{\rho}e_{ij}=\sqrt{\rho}\frac{E_{ij}+E_{ji}}{\sqrt{2}},\sqrt{\rho}e_{ij}^{\alpha}=\sqrt{\rho}\sqrt{-1}\frac{E_{ij}-E_{ji}}{\sqrt{2}}\}$
be the orthonormal basis for $J_{ij}$. We choose the following Cartan basis for
$\mathfrak{su}(n,n)$:

\begin{align*}
   H_{e_{i+\bar{n}_{1}}-e_{j+\bar{n}_{2}}}&=\frac{1}{2}(h_{\eta_{\bar{n}_{1}}e_{ii}+\eta_{\bar{n}_{2}}e_{jj}}-4i[L_{e_{ij}},L_{e_{ij}^{\alpha}}])\\
   H_{e_{i}-e_{n+i}}&=h_{e_{ii}}\\
   E_{e_{i}-e_{n+i}}&=E_{e_{ii}}^{+}=iL_{e_{ii}}-\frac{1}{2}(X_{e_{ii}}-Y_{e_{ii}})
   \\
   E_{-e_{i}+e_{n+i}}&=E_{e_{ii}}^{-}=iL_{e_{ii}}+\frac{1}{2}(X_{e_{ii}}-Y_{e_{ii}})
   \\
   E_{e_{i}-e_{j}}&=\frac{1}{2\sqrt{2}}(h_{e_{ij}+ie_{ij}^{\alpha}}+
   4[L_{e_{ij}+ie_{ij}^{\alpha}},L_{e_{jj}}]) \\
   E_{-e_{i}+e_{j}}&=\frac{1}{2\sqrt{2}}(h_{e_{ij}-ie_{ij}^{\alpha}}-
   4[L_{e_{ij}-ie_{ij}^{\alpha}},L_{e_{jj}}]) \\
   E_{e_{n+i}-e_{n+j}}&=\frac{1}{2\sqrt{2}}(h_{e_{ij}+ie_{ij}^{\alpha}}-
   4[L_{e_{ij}+ie_{ij}^{\alpha}},L_{e_{jj}}]) \\
   E_{-e_{n+i}+e_{n+j}}&=\frac{1}{2\sqrt{2}}(h_{e_{ij}-ie_{ij}^{\alpha}}+
   4[L_{e_{ij}-ie_{ij}^{\alpha}},L_{e_{jj}}]) \\
   E_{e_{i}-e_{n+j}}&=\frac{1}{\sqrt{2}}E_{e_{ij}+ie_{ij}^{\alpha}}^{+}& E_{-e_{i}+e_{n+j}}&=\frac{1}{\sqrt{2}}E_{e_{ij}-ie_{ij}^{\alpha}}^{-}\\
   E_{e_{j}-e_{n+i}}&=\frac{1}{\sqrt{2}}E_{e_{ij}-ie_{ij}^{\alpha}}^{+}& E_{-e_{j}+e_{n+i}}&=\frac{1}{\sqrt{2}}E_{e_{ij}+ie_{ij}^{\alpha}}^{-}\\
 \end{align*}
Where $\eta_{\bar{n}_{1}}=\left\{
                           \begin{array}{ll}
                             1, & \hbox{when $\bar{n}_{1}=0$ ;} \\
                             -1, & \hbox{when $\bar{n}_{1}=n$.}\\
                           \end{array}
                         \right.$

Similarly for $\eta_{\bar{n}_{2}}$.

\hspace{3cm}

Then we get:

 \begin{align*}
      L_{e_{ii}}&=-\frac{i}{2}( E_{e_{i}-e_{n+i}}+ E_{-e_{i}+e_{n+i}})\\
      L_{e_{ij}}&=-\frac{i}{2\sqrt{2}}(E_{e_{i}-e_{n+j}}+E_{-e_{i}+e_{n+j}}+E_{e_{j}-e_{n+i}}+E_{-e_{j}+e_{n+i}})\\
      L_{e_{ij}^{\alpha}}&=-\frac{1}{2\sqrt{2}}(E_{e_{i}-e_{n+j}}-E_{-e_{i}+e_{n+j}}-E_{e_{j}-e_{n+i}}+E_{-e_{j}+e_{n+i}})\\
      X_{e_{ii}}&=\frac{1}{2}(E_{-e_{i}+e_{n+i}}-E_{e_{i}-e_{n+i}}+iH_{e_{i}-e_{n+i}})\\
      Y_{e_{ii}}&=-\frac{1}{2}(E_{-e_{i}+e_{n+i}}-E_{e_{i}-e_{n+i}}-iH_{e_{i}-e_{n+i}})\\
      X_{e_{ij}}&=\frac{1}{2\sqrt{2}}(iE_{e_{i}-e_{j}}+iE_{-e_{i}+e_{j}}+iE_{e_{n+i}-e_{n+j}}+iE_{-e_{n+i}+e_{n+j}}\\
                &+E_{-e_{i}+e_{n+j}}+E_{-e_{j}+e_{n+i}}-E_{e_{i}-e_{n+j}}-E_{e_{j}-e_{n+i}})\\
      Y_{e_{ij}}&=\frac{1}{2\sqrt{2}}(iE_{e_{i}-e_{j}}+iE_{-e_{i}+e_{j}}+iE_{e_{n+i}-e_{n+j}}+iE_{-e_{n+i}+e_{n+j}}\\
               & -E_{-e_{i}+e_{n+j}}-E_{-e_{j}+e_{n+i}}+E_{e_{i}-e_{n+j}}+E_{e_{j}-e_{n+i}})\\
      X_{e_{ij}^{\alpha}}&=\frac{1}{2\sqrt{2}}(E_{e_{i}-e_{j}}+E_{-e_{i}+e_{j}}-E_{e_{n+i}-e_{n+j}}-E_{-e_{n+i}+e_{n+j}}\\
                &+iE_{-e_{i}+e_{n+j}}-iE_{-e_{j}+e_{n+i}}+iE_{e_{i}-e_{n+j}}-iE_{e_{j}-e_{n+i}})\\
      Y_{e_{ij}^{\alpha}}&=\frac{1}{2\sqrt{2}}(E_{e_{i}-e_{j}}+E_{-e_{i}+e_{j}}-E_{e_{n+i}-e_{n+j}}-E_{-e_{n+i}+e_{n+j}}\\
                &-iE_{-e_{i}+e_{n+j}}+iE_{-e_{j}+e_{n+i}}-iE_{e_{i}-e_{n+j}}+iE_{e_{j}-e_{n+i}})\\
  &\text{where}~ 1\leq i<j\leq n.
  \end{align*}

\hspace{2cm}

Let
($\pi$, V) be the unitary  highest weight module of $\mathfrak{co}(J)=\mathfrak{su}(n,n)$ with highest weight $\lambda=(\lambda_{1},\lambda_{2},...,\lambda_{n})$.
Then for any highest weight vector $v \in V$, we have
$\pi (H_{\alpha})v=\tilde{H}_{\alpha}v=\lambda(H_{\alpha})v=\frac{2(\lambda,\alpha)}{(\alpha,\alpha)}v$
and $\pi (E_{\alpha})v=0$ if $\alpha$ is a positive root.

We let  the quadratic relations (Q1) and (Q2) act on any highest weight vector, then we get:

 \begin{equation}\label{6.10}
     (\sum\limits_{1\leq i\leq n}\lambda_{i}-\lambda_{n+i})^{2}+2n(\sum\limits_{1\leq i\leq n}\lambda_{i}-\lambda_{n+i})\\
     =-4a.\end{equation}
 \begin{equation}\label{6.11}
 \sum\limits_{1\leq i\leq n}(\lambda_{i}-\lambda_{n+i})^{2}-(\sum\limits_{1\leq i\leq n}\lambda_{i}-\lambda_{n+i})^{2}
     -2\sum\limits_{1\leq i<j\leq n}(\lambda_{j}-\lambda_{n+i})=0\end{equation}

From the unitarity
we have \begin{equation}\label{612}\lambda_{n}\leq  \lambda_{n-1}\leq  ...\leq  \lambda_{1}\leq
\lambda_{2n}\leq ...\leq \lambda_{n+2}\leq \lambda_{n+1}.\end{equation}

For A-type root system, we have \begin{equation}\label{Atype}
\sum\limits_{1\leq i\leq 2n}\lambda_{i}=0. \end{equation}

%
%
%
%
%
%
%

From the quadratic relations (Q3) and (Q4), we can get:

\begin{eqnarray}\label{6.15}
    \left\{\begin{array}{l}
(\lambda_{i}-\lambda_{n+i})^{2}+(-\sum\limits_{1\leq i\leq n}(\lambda_{i}-\lambda_{n+i})+2-i)(\lambda_{i}-\lambda_{n+i})\\
+(n+2-2i)\lambda_{n+i}
                     -\sum\limits_{1\leq j\leq n}\lambda_{j}
                            +\sum\limits_{1\leq k\leq i-1}(\lambda_{k}+\lambda_{n+k})=0,\\
                       \text{ here we take}~  u=\sqrt{\rho}e_{ii} ~\text{for}~ 1\leq i\leq n \\
\end{array}\right.
\end{eqnarray}

and

\begin{eqnarray}\label{6.16}
    \left\{\begin{array}{l}
(\lambda_{i}-\lambda_{n+i})^{2}+(-\sum\limits_{1\leq i\leq n}(\lambda_{i}-\lambda_{n+i})+2-i)(\lambda_{i}-\lambda_{n+i})\\
+(n+2-2i)\lambda_{n+i}-\sum\limits_{1\leq j\leq n}\lambda_{j}+\sum\limits_{1\leq k\leq i-1}
(\lambda_{k}+\lambda_{n+k})\\
+(\lambda_{j}-\lambda_{n+j})^{2}
    +(-\sum\limits_{1\leq i\leq n}(\lambda_{i}-\lambda_{n+i})+2-j)(\lambda_{j}-\lambda_{n+j})\\
+(n+2-2j)\lambda_{n+j}-\sum\limits_{1\leq j\leq n}\lambda_{j}
    +\sum\limits_{1\leq k\leq
      j-1}(\lambda_{k}+\lambda_{n+k})\\
      +2(\lambda_{i}-\lambda_{n+i})(\lambda_{j}-\lambda_{n+j})    +2(\lambda_{j}-\lambda_{n+i})=0,\\
    \text{ here we take} ~u=e_{ii}+e_{jj} ~\text{for}~ 1\leq i<j\leq n
\end{array}\right.
\end{eqnarray}

From (\ref{6.15}) and (\ref{6.16}) we get:
\begin{equation}\label{6.17}(\lambda_{i}-\lambda_{n+i})(\lambda_{j}-
\lambda_{n+j})=-(\lambda_{j}-\lambda_{n+i}).\end{equation}

If $\lambda_{i_{0}}-\lambda_{n+i_{0}}=0$ for some $i_{0}$, then we
have $\lambda_{j}-\lambda_{n+i_{0}}=0~for ~i_{0}<j\leq n$ and
$\lambda_{i_{0}}-\lambda_{n+k}=0 ~for~1\leq k<i_{0}$ since (\ref{6.17}).
This result and (\ref{612}) imply
$\lambda_{n}=...=\lambda_{1}=\lambda_{2n}=...=\lambda_{n+1}$.
Since the constant  $a$  is not zero, we get a contradiction with (\ref{6.10}).


This implies
$\lambda_{i}-\lambda_{n+i}\leq -1$ for all $i$ since it is an
integer.

%


\hspace{1cm}

Suppose $\lambda_{1}-\lambda_{n+1}\leq -2$ and
$\lambda_{n}-\lambda_{2n}\leq -2$. Then from (\ref{6.17}) we have
$\lambda_{n}-\lambda_{n+1}\leq -2\lambda_{1}+2\lambda_{n+1}$ and
$\lambda_{n}-\lambda_{n+1}\leq
-2\lambda_{n}+2\lambda_{2n}$. $\Rightarrow 2\lambda_{2n}\leq
-\lambda_{n}-\lambda_{n+1}\leq -2\lambda_{1}$. $\Rightarrow
\lambda_{2n}=\lambda_{1}$ since (\ref{612}). Then we must have
$\lambda_{1}-\lambda_{n+1}= -2$, $\lambda_{n}-\lambda_{2n}= -2$, and
$\lambda_{n}-\lambda_{n+1}=-4$.

Then we have
$-2(\lambda_{i}-\lambda_{n+i})=(\lambda_{i}-\lambda_{n+i})(\lambda_{n}-\lambda_{2n})=-\lambda_{n}+\lambda_{n+i}$
and
$-2(\lambda_{i}-\lambda_{n+i})=(\lambda_{1}-\lambda_{n+1})(\lambda_{i}-\lambda_{n+i})=-\lambda_{i}+\lambda_{n+1}$.
We add them together and get
$-4(\lambda_{i}-\lambda_{n+i})=-(\lambda_{n}-\lambda_{n+i}+\lambda_{i}-\lambda_{n+1})=-(4+\lambda_{i}-\lambda_{n+i})$
$\Rightarrow \lambda_{i}-\lambda_{n+i}=-\frac{4}{3}$. This is a
contradiction since all $\lambda_{i}-\lambda_{n+i}$ are integers.

So we may suppose $\lambda_{n}-\lambda_{2n}=-1$. If we also have
$\lambda_{1}-\lambda_{n+1}=-1$, then  $-1
\geq \lambda_{i}-\lambda_{n+i}=\lambda_{i}-\lambda_{n+1}\geq
\lambda_{n}-\lambda_{n+1}=-1$ $\Rightarrow \lambda_{i}-\lambda_{n+i}=-1$ for all $1\leq i \leq n$.  Now we suppose $\lambda_{1}-\lambda_{n+1}=-1-k \leq-2$. Then
$\lambda_{n}-\lambda_{2n}=-1 \Rightarrow
\lambda_{i}-\lambda_{n+i}=\lambda_{n}-\lambda_{n+i}\Rightarrow
\lambda_{i}=\lambda_{n}$ for all $1\leq i \leq n-1$. From
(\ref{6.17}) we get
$(\lambda_{1}-\lambda_{n+1})(\lambda_{j}-\lambda_{n+j})=-\lambda_{j}+\lambda_{n+1}=-\lambda_{1}+\lambda_{n+1}>0$.
$\Rightarrow \lambda_{j}-\lambda_{n+j}=-1$ for all $1<j\leq n$.

The result is similar if we suppose $\lambda_{1}-\lambda_{n+1}=-1$.
\hspace{4cm}

So we can get:

$\lambda=(\underbrace{b,\dots,b}_{n},b+1+k,b+1,...b+1),k=0,1,...$,

or $\lambda= (b,...b,b-k,\underbrace{b+1,...,b+1}_{n}), k=0,1,...$ .

Then from (\ref{Atype}), we have

$\lambda=(\underbrace{-\frac{n+k}{2n},\dots,-\frac{n+k}{2n}}_{n},\frac{n-k}{2n}+k,\frac{n-k}{2n},...,\frac{n-k}{2n}),k=0,1,...$, \\ or $\lambda=(-\frac{n-k}{2n},...,-\frac{n-k}{2n},-\frac{n-k}{2n}-k,\underbrace{\frac{n+k}{2n},...,\frac{n+k}{2n}}_{n}), k=0,1,...$.


\subsection{The case when $J=\mathcal {H}_{n}(\mathbb{H})(n\geq3)$}
\hspace{1cm}

In this case, $\mathfrak{co}(J)=\mathfrak{so}^{\ast}(4n)$. The compact
positive roots are $e_{i}-e_{j}$ with $1\leq i<j\leq 2n$, and the
noncompact roots are $e_{i}+e_{j}$ with $1\leq i<j\leq 2n$. Let
$\{\sqrt{\rho}e_{ij}=\sqrt{\rho}\frac{E_{ij}+E_{ji}}{\sqrt{2}},\sqrt{\rho}e_{ij}^{2}=\sqrt{\rho}\sqrt{-1}\frac{E_{ij}-E_{ji}}{\sqrt{2}},
\sqrt{\rho}e_{ij}^{3}=\sqrt{\rho}j\frac{E_{ij}-E_{ji}}{\sqrt{2}},\sqrt{\rho}e_{ij}^{4}=\sqrt{\rho}k\frac{E_{ij}-E_{ji}}{\sqrt{2}}\}$
be the orthonormal basis for $J_{ij}$. We choose the following Cartan basis for
$\mathfrak{so}^{\ast}(4n)$:

\begin{alignat}{2}\nonumber
    H_{e_{i}-e_{n+i}}&=2i([L_{e_{ij}},L_{e_{ij}^{4}}]-[L_{e_{ij}^{2}},L_{e_{ij}^{3}}])&  H_{e_{i}+e_{n+i}}&=h_{e_{ii}}
   \\\nonumber
   E_{e_{i}+e_{n+i}}&=E_{e_{ii}}^{+}=iL_{e_{ii}}-\frac{1}{2}(X_{e_{ii}}-Y_{e_{ii}}) \\\nonumber E_{-e_{i}-e_{n+i}}&=E_{e_{ii}}^{-}=iL_{e_{ii}}+\frac{1}{2}(X_{e_{ii}}-Y_{e_{ii}}) \\\nonumber
    E_{e_{i}-e_{n+i}}&=[L_{e_{ij}^{2}}+iL_{e_{ij}^{3}},L_{e_{ij}}-iL_{e_{ij}^{4}}]\\\nonumber
   E_{-e_{i}+e_{n+i}}&=[L_{e_{ij}}+iL_{e_{ij}^{4}},L_{e_{ij}^{2}}-iL_{e_{ij}^{3}}]\\\nonumber
   H_{e_{\bar{n}+i}-e_{\bar{n}+j}}&=\frac{1}{2}(h_{e_{ii}-e_{jj}}+4i\eta_{\bar{n}}[L_{e_{ij}},L_{e_{ij}^{4}}])\\\nonumber
   E_{e_{\bar{n}+i}-e_{\bar{n}+j}}&=\frac{1}{2\sqrt{2}}(h_{e_{ij}-i\eta_{\bar{n}}e_{ij}^{4}}+4[L_{e_{ij}-i\eta_{\bar{n}}e_{ij}^{4}},L_{e_{jj}}])\\\nonumber
   E_{-e_{\bar{n}+i}+e_{\bar{n}+j}}&=\frac{1}{2\sqrt{2}}(h_{e_{ij}+i\eta_{\bar{n}}e_{ij}^{4}}-4[L_{e_{ij}+i\eta_{\bar{n}}e_{ij}^{4}},L_{e_{jj}}])\\\nonumber
   H_{e_{\bar{n}+i}-e_{n+j-\bar{n}}}&=\frac{1}{2}(h_{e_{ii}-e_{jj}}-4i\eta_{\bar{n}}[L_{e_{ij}^{2}},L_{e_{ij}^{3}}])\\\nonumber
   E_{e_{\bar{n}+i}-e_{n+j-\bar{n}}}&=\frac{1}{2\sqrt{2}}(h_{e_{ij}^{2}+i\eta_{\bar{n}}e_{ij}^{3}}+4[L_{e_{ij}^{2}+i\eta_{\bar{n}}e_{ij}^{3}},L_{e_{jj}}])\\\nonumber
   E_{-e_{\bar{n}+i}+e_{n+j-\bar{n}}}&=\frac{1}{2\sqrt{2}}(h_{e_{ij}^{2}-i\eta_{\bar{n}}e_{ij}^{3}}-4[L_{e_{ij}^{2}-i\eta_{\bar{n}}e_{ij}^{3}},L_{e_{jj}}])\\\nonumber
   H_{e_{i+\bar{n}}+e_{\bar{n}+j}}&=\frac{1}{2}(h_{e_{ii}+e_{jj}}-4i\eta_{\bar{n}}[L_{e_{ij}^{2}},L_{e_{ij}^{3}}])\\\nonumber
   E_{e_{i+\bar{n}}+e_{\bar{n}+j}}&=\frac{1}{\sqrt{2}}E_{e_{ij}^{2}+i\eta_{\bar{n}}e_{ij}^{3}}^{+}& E_{-e_{i+\bar{n}}-e_{\bar{n}+j}}&=\frac{1}{\sqrt{2}}E_{e_{ij}^{2}-i\eta_{\bar{n}}e_{ij}^{3}}^{-}\\\nonumber
   E_{e_{\bar{n}+i}+e_{n+j-\bar{n}}}&=\frac{1}{\sqrt{2}}E_{e_{ij}-i\eta_{\bar{n}}e_{ij}^{4}}^{+}& E_{-e_{\bar{n}+i}-e_{n+j-\bar{n}}}&=\frac{1}{\sqrt{2}}E_{e_{ij}+i\eta_{\bar{n}}e_{ij}^{4}}^{-}\\\nonumber
    H_{e_{\bar{n}+i}+e_{n+j-\bar{n}}}&=\frac{1}{2}(h_{e_{ii}+e_{jj}}+4i\eta_{\bar{n}}[L_{e_{ij}},L_{e_{ij}^{4}}])\\\nonumber
 \end{alignat}
 Where $\eta_{\bar{n}}=\left\{
                           \begin{array}{ll}
                             1, & \hbox{when $\bar{n}=0$ ;} \\
                             -1, & \hbox{when $\bar{n}=n$.}\\
                           \end{array}
                         \right.$

~\\

Let
($\pi$, V) be the unitary  highest weight module of $\mathfrak{co}(J)=\mathfrak{so}^{\ast}(4n)$ with highest weight $\lambda=(\lambda_{1},\lambda_{2},...,\lambda_{n})$.
Then for any highest weight vector $v \in V$, we have
$\pi (H_{\alpha})v=\tilde{H}_{\alpha}v=\lambda(H_{\alpha})v=\frac{2(\lambda,\alpha)}{(\alpha,\alpha)}v$
and $\pi (E_{\alpha})v=0$ if $\alpha$ is a positive root.

Similar to the previous case, we can write the generators of
$\mathfrak{so}^{\ast}(4n)$ as linear combinations of these Cartan
basis and then from the quadratic relations (Q1) and (Q2), we can get:

\begin{equation}\label{6.18}
(\sum\limits_{1\leq i\leq n}\lambda_{i}+\lambda_{n+i})^{2}+(4n-2)(\sum\limits_{1\leq i\leq n}\lambda_{i}+\lambda_{n+i})=-4a
\end{equation}

and

\begin{equation}\label{6.19}
     \sum\limits_{1\leq i\leq n}(\lambda_{i}+\lambda_{n+i})^{2}-(\sum\limits_{1\leq i\leq n}\lambda_{i}+\lambda_{n+i})^{2}
     -2\sum\limits_{1\leq i<j\leq n}(\lambda_{j}+\lambda_{n+i}+2\lambda_{n+j})=0.\end{equation}

From the
unitarity we have: \begin{equation}\label{6.21}\lambda_{2n}\leq  \lambda_{2n-1}\leq  ...\leq
\lambda_{n+1}\leq
 \lambda_{n}\leq ...\leq -|\lambda_{1}|.\end{equation}

From the quadratic relations (Q3) and (Q4),  we can get:

\begin{eqnarray}\label{6.23}
\left\{\begin{array}{l}
 (\lambda_{i}+\lambda_{n+i})^{2}+(-\sum\limits_{1\leq i\leq n}(\lambda_{i}+\lambda_{n+i})+2-i)(\lambda_{i}+\lambda_{n+i})\\
 -(n-2)\lambda_{n+i}
 -\sum\limits_{1\leq j\leq n}\lambda_{n+j}
 -\sum\limits_{1\leq i\leq n}(\lambda_{i}+\lambda_{n+i})\\
 +\sum\limits_{1\leq k\leq i-1}(\lambda_{k}+\lambda_{n+k})=0,\\
  \text{ here we take} ~u=\sqrt{\rho}e_{ii} ~\text{for}~ 1\leq i\leq n \\
\end{array}\right.
\end{eqnarray}

and

\begin{eqnarray}\label{6.24}
\left\{\begin{array}{l}
(\lambda_{i}+\lambda_{n+i})^{2}+(-\sum\limits_{1\leq i\leq n}(\lambda_{i}+\lambda_{n+i})+2-i)(\lambda_{i}+\lambda_{n+i})\\
-(n-2)\lambda_{n+i}
-\sum\limits_{1\leq l\leq n}\lambda_{n+l}-\sum\limits_{1\leq i\leq n}(\lambda_{i}+\lambda_{n+i})\\
+\sum\limits_{1\leq k\leq i-1}(\lambda_{k}+\lambda_{n+k})\\+
(\lambda_{j}+\lambda_{n+j})^{2}
+(-\sum\limits_{1\leq i\leq n}(\lambda_{i}+\lambda_{n+i})+2-j)(\lambda_{j}+\lambda_{n+j})\\
-(n-2)\lambda_{n+j}
-\sum\limits_{1\leq l\leq n}\lambda_{n+l}
  -\sum\limits_{1\leq i\leq n}(\lambda_{i}+\lambda_{n+i})\\ +\sum\limits_{1\leq k\leq j-1}(\lambda_{k}+\lambda_{n+k})\\
  +2(\lambda_{i}+\lambda_{n+i})(\lambda_{j}+\lambda_{n+j})
  +2(\lambda_{j}+\lambda_{n+i}+2\lambda_{n+j})=0,\\
   \text{ here we take} ~u=e_{ii}+e_{jj} ~\text{for}~ 1\leq i<j\leq n.
\end{array}\right.
\end{eqnarray}

From (\ref{6.23}) and (\ref{6.24}) we get:
\begin{equation}\label{6.25}(\lambda_{i}+\lambda_{n+i})(\lambda_{j}+\lambda_{n+j})\\
=-(\lambda_{j}+\lambda_{n+i}+2\lambda_{n+j})\end{equation}

From
(\ref{6.21}) we know $\lambda_{i}+\lambda_{n+i}\leq 0$. If
$\lambda_{i_{0}}+\lambda_{n+i_{0}}=0$ for some $i_{0}<n$, then
$(\ref{6.25})\Rightarrow 0=\lambda_{n}+\lambda_{n+i_{0}}+2\lambda_{2n}\leq
0\Rightarrow \lambda_{n}=0,\lambda_{n+i_{0}}=0,
\lambda_{i_{0}}=0,~and~ \lambda_{2n}=0 $. Then we get all
$\lambda_{i}=0, \lambda_{n+i}=0$. If $\lambda_{n}+\lambda_{2n}=0$,
then $(\ref{6.21})\Rightarrow $ all $\lambda_{i}=0, \lambda_{n+i}=0$. But
this is a contradiction since the given highest weight module is
nontrivial. So we must have:
$$2\lambda_{n+i}\leq \lambda_{i}+\lambda_{n+i}\leq -1 ~\text{for all}~i.$$

\begin{Lem}
$\lambda_{i}+\lambda_{n+i}=-2 ~for~ all~ i<n.$

\end{Lem}

{\bf Proof.} If $\lambda_{i_{0}}+\lambda_{n+i_{0}}=-1$ for some
$i_{0}<n$, then $(\ref{6.25})\Rightarrow
\lambda_{n}+\lambda_{2n}=\lambda_{n}+\lambda_{n+i_{0}}+2\lambda_{2n}
\Rightarrow 0=\lambda_{n+i_{0}}+\lambda_{2n}\leq
2\lambda_{n+i_{0}}\leq -1$.This contradiction implies all
$\lambda_{i}+\lambda_{n+i}\leq -2$ for $i<n$.

If $\lambda_{i_{0}}+\lambda_{n+i_{0}}\leq -4$ for some $i_{0}<n$,then
$(\ref{6.25})\Rightarrow 4\lambda_{n}+4\lambda_{2n}\geq
\lambda_{n}+\lambda_{n+i_{0}}+2\lambda_{2n} \Rightarrow
\lambda_{n+i_{0}}\leq  3\lambda_{n}+2\lambda_{2n}\leq
2\lambda_{2n}\leq  2\lambda_{n+i_{0}}\leq -1$. This contradiction
implies all $\lambda_{i}+\lambda_{n+i}\geq -3$ for $i<n$.

If $\lambda_{i_{0}}+\lambda_{n+i_{0}}=-3$ for some $i_{0}<n$, then
$(\ref{6.25})\Rightarrow
\lambda_{n}+\lambda_{n+i_{0}}+2\lambda_{2n}=3\lambda_{n}+3\lambda_{2n}\Rightarrow
\lambda_{n+i_{0}}= 2\lambda_{n}+\lambda_{2n}\leq  \lambda_{2n}\leq
\lambda_{n+i_{0}}\leq -\frac{1}{2}\Rightarrow \lambda_{2n}=...=
\lambda_{n+i_{0}}, \lambda_{n}=0\Rightarrow
all~\lambda_{i}=0~since~(\ref{6.21}), and~
\lambda_{2n}=...=\lambda_{n+i_{0}}=-3$. When $i_{0}>1$, we let
$k<i_{0}$, then $(\ref{6.25})\Rightarrow
3(\lambda_{k}+\lambda_{n+k})=\lambda_{i_{0}}+\lambda_{n+k}+2\lambda_{n+i_{0}}
\Rightarrow \lambda_{n+k}=-3$ from the above computations.
 So we always have $\lambda_{2n}=...=\lambda_{n+1}=-3,
\lambda_{n}=...=\lambda_{1}=0$.
But from the quadratic relation (Q4), we have
\begin{eqnarray}\label{4e12}
\left\{\begin{array}{l}
(\lambda_{1}+\lambda_{2}+\lambda_{n+1}+\lambda_{n+2})(n+1+\sum\limits_{1\leq i\leq n}(\lambda_{i}+\lambda_{n+i}))\\
+2\sum\limits_{1\leq i\leq n}(\lambda_{i}+2\lambda_{n+i})-(\lambda_{1}-\lambda_{2}-\lambda_{n+1}+
\lambda_{n+2})^{2}\\
-(n+1)\lambda_{1}-(n-3)\lambda_{2}-5\lambda_{n+1}-\lambda_{n+2}=0,\\
 \text{ here we take} ~u=\sqrt{\rho}e_{12}
\end{array}\right.
\end{eqnarray}
Then we get:
$-6(n+1-3n)+2(-6n)+15+3=0$, i.e. 12=0, a contradiction!
%
%
%

So we must have

 $\lambda_{i}+\lambda_{n+i}=-2$ for all
$1\leq i<n$.   $\boxempty $

\hspace{4cm}

From the above lemma and (\ref{6.25}), we have:\\
$2(\lambda_{i}+\lambda_{n+i})=(\lambda_{i}+\lambda_{n+1}+2\lambda_{n+i})\Rightarrow
\lambda_{i}=\lambda_{n+1}$ for all $1<i\leq n $.

From the above lemma and (\ref{6.21}) we also have
$\lambda_{1}=...=\lambda_{n-1},\lambda_{n+1}=...=\lambda_{2n-1}$.
 So we must have $\lambda_{1}=...=\lambda_{n}=\lambda_{n+1}=...=\lambda_{2n-1}=-1\geq \lambda_{2n}=-1-k.$
i.e. $\lambda=(-1,...,-1,-k-1), k=0,1,...$ .

\hspace{6cm}


\subsection{The case when $J=\mathcal {H}_{3}(\mathbb{O})$}
\hspace{1cm}

In this case, $\mathfrak{co}(J)=\mathfrak{e}_{7(-25)}$. The compact
positive roots are
$\{\pm e_{i}+e_{j} ~with ~1\leq i<j\leq
5\}\bigcup\{\frac{1}{2}(e_{8}-e_{7}-e_{6}+\sum\limits_{1\leq i\leq
5}(-1)^{n(i)}e_{i})|\sum\limits_{1\leq i\leq 5}n(i)~is ~even\}$, and
the noncompact positive roots are $\{\pm e_{i}+e_{6}~ with ~1\leq
i\leq5\}\bigcup\{e_{8}-e_{7}\}\bigcup\{\frac{1}{2}(e_{8}-e_{7}+e_{6}+\sum\limits_{1\leq
i\leq 5}(-1)^{n(i)}e_{i})|\sum\limits_{1\leq i\leq 5}n(i)~is
~odd\}$. Let
$\{\sqrt{\rho}e_{ij}^{\alpha}=\sqrt{\rho}\alpha\frac{E_{ij}-E_{ji}}{\sqrt{2}}|
\alpha=1,i,j,k,l,il,jl,kl \}$ be the orthonormal basis for $J_{ij}$.
We choose the following Cartan
basis for $\mathfrak{e}_{7(-25)}$:

\begin{align*}\nonumber
E_{e_{6}+e_{3}}&=E_{e_{22}}^{+}&    E_{e_{6}-e_{3}}&=E_{e_{11}}^{+}\\\nonumber
E_{e_{6}+e_{1}}&=E_{z_{1}}^{+}
=\frac{1}{\sqrt{2}}E_{e_{12}+\sqrt{-1}e_{12}^{i}}^{+} &    E_{e_{6}+e_{2}}&=E_{z_{2}}^{+}
=\frac{1}{\sqrt{2}}E_{e_{12}^{j}+\sqrt{-1}e_{12}^{k}}^{+}\\\nonumber
E_{e_{6}+e_{4}}&=E_{z_{4}}^{+}
=\frac{1}{\sqrt{2}}E_{e_{12}^{l}+\sqrt{-1}e_{12}^{il}}^{+}&    E_{e_{6}+e_{5}}&=E_{z_{5}}^{+}
=\frac{1}{\sqrt{2}}E_{e_{12}^{jl}-\sqrt{-1}e_{12}^{kl}}^{+}\\\nonumber
E_{e_{6}-e_{i}}&=E_{\bar{z_{i}}}^{+}~ \text{for}~i=1,2,4,5\\\nonumber
\end{align*}
\begin{align*}\nonumber
E_{e_{i}+e_{j}}&=[L_{z_{i}},L_{z_{j}}]  ~\text{for}~1\leq i<j\leq 5~\text{and}~i\neq 3, j\neq 3\\\nonumber
\end{align*}
\begin{align*}\nonumber
E_{-e_{i}+e_{j}}&=[L_{\bar{z_{i}}},L_{z_{j}}]& E_{-e_{i}-e_{j}}&=[L_{\bar{z_{i}}},L_{\bar{z_{j}}}]\\\nonumber
E_{e_{i}-e_{j}}&=[L_{z_{i}},L_{\bar{z_{j}}}]& E_{e_{i}\pm e_{3}}&=\frac{1}{2}(h_{z_{i}}\pm 4[L_{z_{i}},L_{e_{11}}])\\\nonumber
E_{-e_{i}\pm e_{3}}&=\frac{1}{2}(h_{\bar{z_{i}}}\pm 4[L_{\bar{z_{i}}},L_{e_{11}}])\\\nonumber
\end{align*}
\begin{align*}\nonumber
E_{-\alpha}&=E_{\bar{z}}^{-}~ ~\text {if}~E_{\alpha}=E_{z}^{+} ~\text{for noncompact positive  root}~ \alpha \\\nonumber
\end{align*}
\begin{align*}\nonumber
E_{\beta_{(-5)}}&=E_{\frac{1}{2}(e_{8}-e_{7}
+e_{6}+e_{3}-e_{5}+e_{4}+e_{2}+e_{1})}
=\frac{1}{\sqrt{2}}E_{e_{23}^{jl}+\sqrt{-1}e_{23}^{kl}}^{+}\\\nonumber
E_{\beta_{(+5)}}&=E_{\frac{1}{2}(e_{8}-e_{7}+e_{6}
+e_{3}+e_{5}-e_{4}-e_{2}-e_{1})}
=\frac{1}{\sqrt{2}}E_{e_{23}^{jl}-\sqrt{-1}e_{23}^{kl}}^{+}\\\nonumber
E_{\beta_{(-4)}}&=E_{\frac{1}{2}(e_{8}-e_{7}+e_{6}
+e_{3}+e_{5}-e_{4}+e_{2}+e_{1})}=\frac{1}{\sqrt{2}}E_{e_{23}^{l}
-\sqrt{-1}e_{23}^{il}}^{+}\\\nonumber
E_{\beta_{(+4)}}&=E_{\frac{1}{2}(e_{8}-e_{7}
+e_{6}+e_{3}-e_{5}+e_{4}-e_{2}-e_{1})}
=\frac{1}{\sqrt{2}}E_{e_{23}^{l}+\sqrt{-1}e_{23}^{il}}^{+}\\\nonumber
\end{align*}
\begin{align*}\nonumber
E_{\beta_{(-2)}}&=\frac{1}{\sqrt{2}}E_{e_{23}^{j}
-\sqrt{-1}e_{23}^{k}}^{+}& E_{\beta_{(+2)}}&=\frac{1}{\sqrt{2}}E_{e_{23}^{j}
+\sqrt{-1}e_{23}^{k}}^{+}\\\nonumber
E_{\beta_{(-1)}}&=\frac{1}{\sqrt{2}}E_{e_{23}
+\sqrt{-1}e_{23}^{i}}^{+}& E_{\beta_{(+1)}}&=\frac{1}{\sqrt{2}}E_{e_{23}
-\sqrt{-1}e_{23}^{i}}^{+}\\\nonumber
\end{align*}
\begin{align*}\nonumber
E_{\mu_{(4-)}}&=\frac{1}{\sqrt{2}}E_{\frac{1}{2}(e_{8}-e_{7}
+e_{6}-e_{3}-e_{5}-e_{4}-e_{2}-e_{1})}
=\frac{1}{\sqrt{2}}E_{e_{13}-\sqrt{-1}e_{13}^{i}}^{+}\\\nonumber
E_{\mu_{(4+)}}&=\frac{1}{\sqrt{2}}E_{\frac{1}{2}(e_{8}
-e_{7}+e_{6}-e_{3}+e_{5}+e_{4}+e_{2}+e_{1})}
=\frac{1}{\sqrt{2}}E_{e_{13}+\sqrt{-1}e_{13}^{i}}^{+}\\\nonumber
E_{\mu_{(+1+2)}}&=\frac{1}{\sqrt{2}}E_{\frac{1}{2}(e_{8}
-e_{7}+e_{6}-e_{3}-e_{5}-e_{4}+e_{2}+e_{1})}
=\frac{1}{\sqrt{2}}E_{e_{13}^{j}+\sqrt{-1}e_{13}^{k}}^{+}\\\nonumber
E_{\mu_{(-1-2)}}&=\frac{1}{\sqrt{2}}E_{\frac{1}{2}(e_{8}-e_{7}
+e_{6}-e_{3}+e_{5}+e_{4}-e_{2}-e_{1})}
=\frac{1}{\sqrt{2}}E_{e_{13}^{j}-\sqrt{-1}e_{13}^{k}}^{+}\\\nonumber
E_{\mu_{(+1+4)}}&=\frac{1}{\sqrt{2}}E_{\frac{1}{2}(e_{8}-e_{7}
+e_{6}-e_{3}-e_{5}+e_{4}-e_{2}+e_{1})}
=\frac{1}{\sqrt{2}}E_{e_{13}^{l}+\sqrt{-1}e_{13}^{il}}^{+}\\\nonumber
E_{\mu_{(-1-4)}}&=\frac{1}{\sqrt{2}}E_{\frac{1}{2}(e_{8}-e_{7}
+e_{6}-e_{3}+e_{5}-e_{4}+e_{2}-e_{1})}
=\frac{1}{\sqrt{2}}E_{e_{13}^{l}-\sqrt{-1}e_{13}^{il}}^{+}\\\nonumber
E_{\mu_{(+2+4)}}&=\frac{1}{\sqrt{2}}E_{\frac{1}{2}(e_{8}-e_{7}
+e_{6}-e_{3}-e_{5}+e_{4}+e_{2}-e_{1})}
=\frac{1}{\sqrt{2}}E_{e_{13}^{jl}+\sqrt{-1}e_{13}^{kl}}^{+}\\\nonumber
E_{\mu_{(-2-4)}}&=\frac{1}{\sqrt{2}}E_{\frac{1}{2}(e_{8}-e_{7}
+e_{6}-e_{3}+e_{5}-e_{4}-e_{2}+e_{1})}
=\frac{1}{\sqrt{2}}E_{e_{13}^{jl}-\sqrt{-1}e_{13}^{kl}}^{+}\\\nonumber
H_{\alpha}&=[E_{\alpha},E_{-\alpha}] ~\text{if $\alpha$  is  compact}\\\nonumber
H_{\alpha}&=-[E_{\alpha},E_{-\alpha}] ~~\text{if $\alpha$  is  noncompact}\nonumber
\end{align*}

\hspace{3cm}

Let
($\pi$, V) be the unitary  highest weight module of $\mathfrak{co}(J)=\mathfrak{e}_{7(-25)}$ with highest weight $\lambda=(\lambda_{1},\lambda_{2},...,\lambda_{n})$.
Then for any highest weight vector $v \in V$, we have
$\pi (H_{\alpha})v=\tilde{H}_{\alpha}v=\lambda(H_{\alpha})v=\frac{2(\lambda,\alpha)}{(\alpha,\alpha)}v$
and $\pi (E_{\alpha})v=0$ if $\alpha$ is a positive root.

Similar to the previous case, we can write the generators of
$\mathfrak{e}_{7(-25)}$ as linear combinations of these Cartan basis
and then from the quadratic relations (Q1) and (Q2), we can get:

\begin{equation}\label{6.26}
      (2\lambda_{6}+\lambda_{8}-\lambda_{7})^{2}+18(2\lambda_{6}+\lambda_{8}-\lambda_{7})+4a=0\\
\end{equation}
and

\begin{eqnarray}\label{6.27}
\left\{\begin{array}{l}(2\lambda_{6}+\lambda_{8}-\lambda_{7})^{2}-(\lambda_{6}-\lambda_{3})^{2} -(\lambda_{6}+\lambda_{3})^{2}\\-
   (\lambda_{8}-\lambda_{7})^{2}
+(24\lambda_{6}-4\lambda_{3}-2\lambda_{4}-2\lambda_{5})=0
\end{array}\right.
\end{eqnarray}

%
%
%

By unitarity we conclude  that

\begin{equation}\label{u4}
-\lambda_{6}\geq  \lambda_{5}\geq  \lambda_{4}\geq  \lambda_{3}\geq
\lambda_{2}\geq |\lambda_{1}|\geq 0, \lambda_{8}\leq 0
\leq \lambda_{7}.\end{equation}

By unitarity, we also have

$\lambda_{8}-\lambda_{7}-\sum\limits_{1\leq j\leq
5,j\neq i}\lambda_{j}\geq (\lambda_{6}-\lambda_{i})$.

Take sum on $i$, then we have:

$10\lambda_{7}+5\lambda_{6}+3\sum\limits_{1\leq j\leq
5}\lambda_{j}\leq 0.$ \\

$\Rightarrow10\lambda_{7}+5\lambda_{6}\leq
-3\sum\limits_{1\leq j\leq 5}\lambda_{j}\leq 0$ \\

$\Rightarrow 2\lambda_{7}+\lambda_{6}\leq 0$.

 If $\lambda_{6}=0$, then $0\leq-\lambda_{8}=\lambda_{7}\leq -\frac{\lambda_{6}}{2}=0$.   From
(\ref{u4}), we can get:

$0=\lambda_{1}=\lambda_{2}=\lambda_{3}=\lambda_{4}=\lambda_{5}$.

This implies:
$\lambda=(0,0,0,0,0,0,0,0)$.

So we can assume $\lambda_{6}<0$.

From the quadratic relation (Q3), we have
\begin{equation}
\lambda_{6}(4-2\lambda_{7})=0, \text{ here we take} ~u=\sqrt{\rho}e_{33}.
\end{equation}
This implies $\lambda_{7}=2, \lambda_{8}=-2$.

From the quadratic relation (Q4), we have
\begin{equation}
\lambda_{7}(4+\lambda_{6})=0, \text{ here we take} ~u=\sqrt{\rho}e_{33}.
\end{equation}
This implies $\lambda_{6}=-4$.

Then from (\ref{6.27}) and (\ref{u4}), we have
$0=\lambda_{1}=\lambda_{2}=\lambda_{3}=\lambda_{4}=\lambda_{5}$.

\hspace{3cm}

%
%
%
%
%

So $\lambda=(0,0,0,0,0,-4,2,-2)$.

\hspace{6cm}

\section{Proof of the main theorem-the hermitian cases (Part II)}
In this section, we will show that for any  unitary highest weight $\mathfrak{co}(J)$-module $L(\lambda)$, if it has the smallest positive Gelfand-Kirillov dimension, then it must satisfy the quadratic relation $(\mathrm{Q}1)$ for some value $a$ only depending on $\lambda=\lambda(k)$ in Corollary \ref{weight k}.

We have the following lemma.
\begin{Lem}
Our quadratic relation $(Q1)$ is equivalent to the following relation:
$$\sum\limits_{1\leq \alpha \leq D}\{ \tilde{X}_{e_{\alpha}},\tilde{Y}_{e_{\alpha}}\}=2\rho(\tilde{L}_{e}^{2}+a). $$

We denote it by $(Q1)'$.
\end{Lem}

{\bf Proof.}
We suppose the relation $(Q1)$ holds, then we compute $[(Q1), \tilde{X}_{e}]$ and get (\ref{5.1}):
$$
\sum\limits_{1\leq \alpha \leq D}\{ \tilde{X}_{e_{\alpha}},\tilde{L}_{e_{\alpha}}\}=\rho\{ \tilde{X}_{e},\tilde{L}_{e}\}.$$

Then from $[(\ref{5.1}), \tilde{Y}_{e}]$ we can get:
$$4\sum\limits_{1\leq\alpha\leq D}\tilde{L}_{e_{\alpha}}^{2}+\sum\limits_{1\leq\alpha\leq D}\{ \tilde{X}_{e_{\alpha}},\tilde{Y}_{e_{\alpha}}\}=4\rho\tilde{L}_{e}^{2}+\rho\{X_{e}, Y_{e}\},$$

which implies the relation $(Q1)'$.

Conversely, if the relation $(Q1)'$ holds, the we compute $[(Q1)', \tilde{X}_{e}]$ and get (\ref{5.1}):
$$
\sum\limits_{1\leq \alpha \leq D}\{ \tilde{X}_{e_{\alpha}},\tilde{L}_{e_{\alpha}}\}=\rho\{ \tilde{X}_{e},\tilde{L}_{e}\}.$$

Then from $[(\ref{5.1}), \tilde{Y}_{e}]$ , we find that the relation $(Q1)$ holds. $\boxempty $

\vspace{1cm}

If a (non-trivial) unitary highest weight $\mathfrak{co}(J)$-module has the smallest positive Gelfand-Kirillov dimension,
then from the previous lemma we need to show that
$\big(\sum\limits_{1\leq \alpha \leq D}\{ \tilde{X}_{e_{\alpha}},\tilde{Y}_{e_{\alpha}}\}-2\rho(\tilde{L}_{e}^{2}+a)\big)v=0$ for any element $v\in V$.

If we denote $Q_{1}'\triangleq \big( \sum\limits_{1\leq \alpha \leq D}\{ \tilde{X}_{e_{\alpha}},\tilde{Y}_{e_{\alpha}}\}-2\rho(\tilde{L}_{e}^{2}+a)\big)$, then from the above case-by-case computation we know:
$ Q_{1}'v_{\lambda}=0$ for any highest weight vector $v_{\lambda}~ in~ V$.


From the TKK commutation relations and the property of highest weight module, we know any element $v\in V$ is a linear combination of  the form $ (\tilde{L}_{e_{\alpha_{1}}}\cdot\cdot\cdot \tilde{L}_{e_{\alpha_{k}}}\tilde{X}_{e_{\beta_{1}}}\cdot\cdot\cdot \tilde{X}_{e_{\beta_{p}}}   \tilde{Y}_{e_{\gamma_{1}}}\cdot\cdot\cdot \tilde{Y}_{e_{\gamma_{q}}})v_{\lambda}$, with $k, p, q\geq 0$. So we need to show \begin{equation} \label{Q1'}
Q_{1}' (\tilde{L}_{e_{\alpha_{1}}}\cdot\cdot\cdot \tilde{L}_{e_{\alpha_{k}}}\tilde{X}_{e_{\beta_{1}}}\cdot\cdot\cdot \tilde{X}_{e_{\beta_{p}}}   \tilde{Y}_{e_{\gamma_{1}}}\cdot\cdot\cdot \tilde{Y}_{e_{\gamma_{q}}})v_{\lambda}=0,
\end{equation}
or equivalently,
\begin{equation}\label{Q1''}
[[[[[[[[\cdot \cdot\cdot[Q_{1}',\tilde{L}_{e_{\alpha_{1}}}],\cdot\cdot\cdot], \tilde{L}_{e_{\alpha_{k}}}],\tilde{X}_{e_{\beta_{1}}}],\cdot\cdot\cdot ], \tilde{X}_{e_{\beta_{p}}}],   \tilde{Y}_{e_{\gamma_{1}}}],\cdot\cdot\cdot], \tilde{Y}_{e_{\gamma_{q}}}]v_{\lambda}=0.
\end{equation}


Let $J$ be a hermitian type simple Euclidean Jordan algebra $\mathcal{H}_{n}(K)$, then we can take the following orthonormal basis:\\
 $\{\sqrt{\rho}e_{pp}, \sqrt{\rho}e_{ab}^{\mu}| 1\leq p\leq \rho, 1\leq a<b\leq \rho, e_{pp}=E_{pp}, e_{ab}=\frac{1}{\sqrt{2}}(E_{ab}+E_{ba}), e_{ab}^{\mu}=\frac{a_{\mu}}{\sqrt{2}} (E_{ab}-E_{ba}) ~for~ 2\leq\mu \leq d~ and~ a_{2}=i, a_{3}=j, a_{4}=k, a_{5}=l, a_{6}=il, a_{7}=jl, a_{8}=kl\},$ where $E_{ab}=(x_{pq})_{n\times n}, x_{pq}=\delta_{ap} \delta_{bq}.$

  We denote this basis by $O_{J}$.

Then we have the following lemma.

\begin{Lem}
For any $e_{\beta}, e_{\beta_{1}},e_{\gamma},...\in O_{J}$, we simply write $\tilde{L}_{e_{\beta}}$ by $\tilde{L}_{\beta}$, $\tilde{X}_{e_{\beta}}$ by $\tilde{X}_{\beta}$, $\tilde{Y}_{e_{\beta}}$ by $\tilde{Y}_{\beta}$, then

\begin{enumerate}
  \item $[Q_{1}', \tilde{L}_{\beta}]=0;$
  \item $A_{1}:=[Q_{1}', \tilde{X}_{\beta}]=2\sum\limits_{1\leq \alpha \leq D}\{\tilde{X}_{\alpha}, \tilde{S}_{\beta \alpha}\}-2\rho\{\tilde{L}_{e}, \tilde{X}_{\beta} \};$
  \item $A_{2}:=[Q_{1}', \tilde{Y}_{\gamma}]=-2\sum\limits_{1\leq \alpha \leq D}\{\tilde{Y}_{\alpha},\tilde{ S}_{\alpha\gamma}\}+2\rho\{\tilde{L}_{e}, \tilde{Y}_{\gamma}\};$
   \item $A_{11}:=[[Q_{1}', \tilde{X}_{\beta}], \tilde{X}_{\beta_{1}}]=4\sum\limits_{1\leq \alpha \leq D}\tilde{X}_{\alpha}\tilde{X}_{S_{\beta\alpha}(\beta_{1})}-4\rho \tilde{X}_{\beta}\tilde{X}_{\beta_{1}},
       ~~[A_{11}, \tilde{X}_{\beta_{2}}]=0;$

           \hspace{1cm}

 \item $A_{22}:=[[Q_{1}', \tilde{Y}_{\gamma}], \tilde{Y}_{\gamma_{1}}]=4\sum\limits_{1\leq \alpha \leq D}\tilde{Y}_{\alpha}\tilde{Y}_{S_{\gamma\alpha}(\gamma_{1})}-4\rho \tilde{Y}_{\gamma}\tilde{Y}_{e_{\gamma_{1}}},
      ~~[A_{22}, \tilde{Y}_{\gamma_{2}}]=0;$

 \item \begin{align*} A_{12}:=&[[Q_{1}', \tilde{X}_{\beta}], \tilde{Y}_{\gamma}]\\
 =&-2\sum\limits_{1\leq \alpha \leq D}\{\tilde{X}_{\alpha}, \tilde{Y}_{S_{\beta\alpha}(\gamma)}\}-4\sum\limits_{1\leq \alpha \leq D}\{\tilde{S}_{\alpha \gamma}, \tilde{S}_{\beta\alpha}\}\\
     &+2\rho\{\tilde{X}_{\beta},\tilde{Y}_{\gamma}\}+4\rho\{\tilde{S}_{\beta \gamma}, \tilde{L}_{e}\};
     \end{align*}

\item \begin{align*} A_{122}:=&[[[Q_{1}', \tilde{X}_{\beta}], \tilde{Y}_{\gamma}],\tilde{Y}_{\gamma_{1}}]\\
=&4\sum\limits_{1\leq \alpha \leq D}\{\tilde{S}_{\alpha \gamma}, \tilde{Y}_{S_{\alpha\beta}(\gamma_{1})}\}+4\sum\limits_{1\leq \alpha \leq D}\{\tilde{S}_{\alpha \gamma_{1}}, \tilde{Y}_{S_{\alpha\beta}(\gamma)}\}\\
    &+4\sum\limits_{1\leq \alpha \leq D}\{\tilde{S}_{\beta\alpha}, \tilde{Y}_{S_{\gamma\alpha}(\gamma_{1})}\}-4\rho\{\tilde{S}_{\beta \gamma_{1}}, \tilde{Y}_{\gamma}\}\\
    &-4\rho\{\tilde{S}_{\beta \gamma}, \tilde{Y}_{\gamma_{1}}\}-4\rho\{\tilde{L}_{e}, \tilde{Y}_{S_{\gamma\beta}(\gamma_{1})}\};
\end{align*}
\item \begin{align*}A_{1222}:=&[[[[Q_{1}', \tilde{X}_{\beta}], \tilde{Y}_{\gamma}],\tilde{Y}_{\gamma_{1}}],\tilde{Y}_{\gamma_{2}}]\\
=&-8\sum\limits_{1\leq \alpha \leq D}\tilde{Y}_{S_{\gamma\alpha}(\gamma_{2})}\tilde{Y}_{S_{\alpha\beta}(\gamma_{1})}-8\sum\limits_{1\leq \alpha \leq D}\tilde{Y}_{S_{\gamma_{1}\alpha(\gamma_{2})}}\tilde{Y}_{S_{\alpha\beta}(\gamma)}\\
    &-8\sum\limits_{1\leq \alpha \leq D}\tilde{Y}_{S_{\alpha\beta}(\gamma_{2})}\tilde{Y}_{S_{\gamma\alpha}(\gamma_{1})}+8\rho \tilde{Y}_{S_{\gamma_{1}\beta}(\gamma_{2})} \tilde{Y}_{\gamma}\\
    &+8\rho \tilde{Y}_{S_{\gamma\beta}(\gamma_{2})}\tilde{Y}_{\gamma_{1}}+8\rho \tilde{Y}_{\gamma_{2}}\tilde{Y}_{S_{\gamma\beta}(\gamma_{1})}\\
    ~~~~~[A_{1222}, \tilde{Y}_{\gamma_{3}}]=0;
\end{align*}
\item \begin{align*} A_{112}:=&[A_{11}, \tilde{Y}_{\gamma}]\\
=&-4\sum\limits_{1\leq \alpha \leq D}\{\tilde{X}_{\alpha},\tilde{S}_{S_{\beta \alpha}(\beta_{1})\cdot\gamma}\}-4\sum\limits_{1\leq \alpha \leq D}\{\tilde{X}_{S_{\beta \alpha}(\beta_{1})},\tilde{S}_{\alpha \gamma}\}\\
    &+4\rho\{\tilde{X}_{\beta},\tilde{S}_{\beta_{1}\gamma}\}+4\rho\{\tilde{X}_{\beta_{1}},\tilde{S}_{\beta \gamma}\};
\end{align*}
\item We simply write $S_{uv}(z)$ by $\{uvz\}$, then
\begin{align*} A_{1122}:=&[A_{112}, \tilde{Y}_{\gamma_{1}}]\\
=&4\sum\limits_{1\leq \alpha \leq D}\{\tilde{X}_{\alpha},\tilde{Y}_{\{\gamma \{\beta \alpha \beta_{1}\}\gamma_{1}\}}\}+8\sum\limits_{1\leq \alpha \leq D}\{\tilde{S}_{\alpha \gamma_{1}},\tilde{S}_{\{\beta \alpha \beta_{1}\}\gamma}\}\\
    &+4\sum\limits_{1\leq \alpha \leq D}\{\tilde{X}_{S_{\beta \alpha}(\beta_{1})},\tilde{Y}_{S_{ \gamma \alpha}(\gamma_{1})}\}+8\sum\limits_{1\leq \alpha \leq D}\{\tilde{S}_{\{\beta \alpha \beta_{1}\} \gamma_{1}},\tilde{S}_{\alpha \gamma}\}\\
    &-4\rho\{\tilde{X}_{\beta},\tilde{Y}_{\{\gamma\beta_{1}\gamma_{1}\}}\}-8\rho\{\tilde{S}_{\beta\gamma_{1}},\tilde{S}_{\beta_{1}\gamma}\}\\
    &-4\rho\{\tilde{X}_{\beta_{1}},\tilde{Y}_{\{\gamma\beta \gamma_{1}\}}\}-8\rho\{\tilde{S}_{\beta_{1}\gamma_{1}},\tilde{S}_{\beta \gamma}\};
\end{align*}
\item \begin{align*}A_{11222}:=&[A_{1122}, \tilde{Y}_{\gamma_{2}}]\\
=&-8\sum\limits_{1\leq \alpha \leq D}\{\tilde{S}_{\alpha\gamma_{2}},\tilde{Y}_{\{\gamma \{\beta \alpha \beta_{1}\}\gamma_{1}\}}\}-8\sum\limits_{1\leq \alpha \leq D}\{\tilde{S}_{\alpha \gamma_{1}},\tilde{Y}_{\{\gamma\{\beta \alpha \beta_{1}\}\gamma_{2}\}}\}\\
&-8\sum\limits_{1\leq \alpha \leq D}\{\tilde{S}_{\{\beta \alpha \beta_{1}\}\gamma},\tilde{Y}_{\{\gamma_{1}\alpha \gamma_{2}\}}\}
    -8\sum\limits_{1\leq \alpha \leq D}\{\tilde{S}_{\{\beta \alpha \beta_{1}\}\gamma_{2}},\tilde{Y}_{S_{ \gamma \alpha}(\gamma_{1})}\}\\
    &-8\sum\limits_{1\leq \alpha \leq D}\{\tilde{S}_{\{\beta \alpha \beta_{1}\} \gamma_{1}},\tilde{Y}_{\{\gamma \alpha \gamma_{2}\}}\}-8\sum\limits_{1\leq \alpha \leq D}\{\tilde{S}_{\alpha \gamma},\tilde{Y}_{\{\gamma_{1}\{\beta \alpha \beta_{1}\} \gamma_{2}\}}\}\\
    &+8\rho\{\tilde{S}_{\beta\gamma_{2}},\tilde{Y}_{\{\gamma\beta_{1}\gamma_{1}\}}\}+8\rho\{\tilde{S}_{\beta\gamma_{1}},\tilde{Y}_{\{\gamma\beta_{1}\gamma_{2}\}}\}\\
    &+8\rho\{\tilde{S}_{\beta_{1}\gamma},\tilde{Y}_{\{\gamma_{1}\beta\gamma_{2}\}}\}
    +8\rho\{\tilde{S}_{\beta_{1}\gamma_{2}},\tilde{Y}_{\{\gamma\beta \gamma_{1}\}}\}\\
    &+8\rho\{\tilde{S}_{\beta_{1}\gamma_{1}},\tilde{Y}_{\{\gamma\beta \gamma_{2}\}}\}+8\rho\{\tilde{S}_{\beta \gamma},\tilde{Y}_{\{\gamma_{1}\beta_{1}\gamma_{2}\}}\};
\end{align*}
\item\begin{align*}A_{112222}:=&[A_{11222}, \tilde{Y}_{\gamma_{3}}]\\
=&16\sum\limits_{1\leq \alpha \leq D}\tilde{Y}_{\{\gamma_{2}\alpha\gamma_{3}\}}\tilde{Y}_{\{\gamma \{\beta \alpha \beta_{1}\}\gamma_{1}\}}+16\sum\limits_{1\leq \alpha \leq D}\tilde{Y}_{\{\gamma_{1}\alpha \gamma_{3}\}}\tilde{Y}_{\{\gamma\{\beta \alpha \beta_{1}\}\gamma_{2}\}}\\
&+16\sum\limits_{1\leq \alpha \leq D}\tilde{Y}_{\{\gamma\{\beta \alpha \beta_{1}\}\gamma_{3}\}}\tilde{Y}_{\{\gamma_{1}\alpha \gamma_{2}\}}
    +16\sum\limits_{1\leq \alpha \leq D}\{\tilde{Y}_{\{\gamma_{2}\{\beta \alpha \beta_{1}\}\gamma_{3}\}}\tilde{Y}_{S_{ \gamma \alpha}(\gamma_{1})}\\
    &+16\sum\limits_{1\leq \alpha \leq D}\{\tilde{Y}_{\{\gamma_{1}\{\beta \alpha \beta_{1}\} \gamma_{3}\}}\tilde{Y}_{\{\gamma \alpha \gamma_{2}\}}+16\sum\limits_{1\leq \alpha \leq D}\tilde{Y}_{\{\gamma\alpha \gamma_{3}\}}\tilde{Y}_{\{\gamma_{1}\{\beta \alpha \beta_{1}\} \gamma_{2}\}}\\
    &-16\rho\tilde{Y}_{\{\gamma_{2}\beta\gamma_{3}\}}\tilde{Y}_{\{\gamma\beta_{1}\gamma_{1}\}}-16\rho\tilde{Y}_{\{\gamma_{1}\beta\gamma_{3}\}}\tilde{Y}_{\{\gamma\beta_{1}\gamma_{2}\}}\\
    &-16\rho\tilde{Y}_{\{\gamma\beta_{1}\gamma_{3}\}}\tilde{Y}_{\{\gamma_{1}\beta\gamma_{2}\}}
    -16\rho\tilde{Y}_{\{\gamma_{2}\beta_{1}\gamma_{3}\}}\tilde{Y}_{\{\gamma\beta \gamma_{1}\}}\\
    &-16\rho\tilde{Y}_{\{\gamma_{1}\beta_{1}\gamma_{3}\}}\tilde{Y}_{\{\gamma\beta \gamma_{2}\}}-16\rho\tilde{Y}_{\{\gamma\beta \gamma_{3}\}}\tilde{Y}_{\{\gamma_{1}\beta_{1}\gamma_{2}\}};\\
[A_{112222}, \tilde{Y}_{\gamma_{4}}]=0;
\end{align*}
\item \begin{align*}
A_{1}v_{\lambda}&=0,&   A_{2}v_{\lambda}&=0;\\
A_{11}v_{\lambda}&=0,&  A_{22}v_{\lambda}&=0;\\
A_{12}v_{\lambda}&=0, &  A_{122}v_{\lambda}&=0;\\
A_{1222}v_{\lambda}&=0, & A_{112}v_{\lambda}&=0;\\
A_{1122}v_{\lambda}&=0, & A_{11222}v_{\lambda}&=0;\\
A_{112222}v_{\lambda}&=0;
\end{align*}
\end{enumerate}
\end{Lem}

\begin{Rem}
For each hermitian type simple Euclidean Jordan algebra $\mathcal {H}_{n}(\mathbb{K})$, its basis $O_{J}$ is a finite set. So the number of the   operators in the form $A_{\alpha}$ in the above lemma must be finite.
\end{Rem}

The proof for this lemma is a straightforward computation, so we skip it.

\hspace{2cm}

{\bf Proof of equation \ref{Q1''}}.  We use the same notation $A_{\alpha}$  with some operator in the lemma to mean they are in the same form.  From the above lemma, we only need to show
\begin{equation}\label{Q1'''}
[[[[[\cdot \cdot\cdot[Q_{1}',\tilde{X}_{\beta_{1}}],\cdot\cdot\cdot ], \tilde{X}_{\beta_{p}}],   \tilde{Y}_{\gamma_{1}}],\cdot\cdot\cdot], \tilde{Y}_{\gamma_{q}}]v_{\lambda}=0,\end{equation}
 since $[Q_{1}', \tilde{L}_{\beta}]=0$.

Also we know $[[[[Q_{1}',\tilde{X}_{\beta_{1}}], \tilde{X}_{\beta_{2}}],\tilde{X}_{\beta_{3}}]=[A_{11},\tilde{X}_{\beta_{3}}]=0.$

So we only need to show the cases  $p=0$, $p=1$ and $p=2$.

For $p=0$, we have
\begin{equation}[[[Q_{1}',\tilde{Y}_{\gamma_{1}}],\cdot\cdot\cdot], \tilde{Y}_{\gamma_{q}}]v_{\lambda}
=\left\{
   \begin{array}{ll}
      A_{2}v_{\lambda}, & {if~ q=1;} \\
      A_{22}v_{\lambda}, & {if~ q=2;} \\
      0, & {if~ q\geq 3.}
   \end{array}
 \right\}
=0.\end{equation}

For $p=1$, we have
\begin{align*}
&[[[[Q_{1}',\tilde{X}_{\beta_{1}}],\tilde{Y}_{\gamma_{1}}],\cdot\cdot\cdot], \tilde{Y}_{\gamma_{q}}]v_{\lambda}\\
=&[[[A_{1},\tilde{Y}_{\gamma_{1}}],\cdot\cdot\cdot], \tilde{Y}_{\gamma_{q}}]v_{\lambda}\\
=&\left\{
   \begin{array}{ll}
      A_{1}v_{\lambda}, & {if~ q=0;}\\
      A_{12}v_{\lambda}, & {if~ q=1;} \\
      A_{122}v_{\lambda}, & {if~ q=2;} \\
      A_{1222}v_{\lambda}, & {if~ q=3;}\\
         0, & {if~ q\geq 4.}
   \end{array}
 \right\}\\
=&0.\end{align*}

For $p=2$, we have
\begin{align*}
&[[[[[Q_{1}',\tilde{X}_{\beta_{1}}],\tilde{X}_{\beta_{2}}],\tilde{Y}_{\gamma_{1}}],\cdot\cdot\cdot], \tilde{Y}_{\gamma_{q}}]v_{\lambda}\\
=&[[[A_{11},\tilde{Y}_{\gamma_{1}}],\cdot\cdot\cdot], \tilde{Y}_{\gamma_{q}}]v_{\lambda}\\
=&\left\{
   \begin{array}{ll}
      A_{11}v_{\lambda}, & {if~ q=0;}\\
      A_{112}v_{\lambda}, & {if~ q=1;} \\
      A_{1122}v_{\lambda}, & {if~ q=2;} \\
      A_{11222}v_{\lambda}, & {if~ q=3;}\\
      A_{112222}v_{\lambda}, & {if~ q=4;}\\
     0, & {if~ q\geq 5.}
   \end{array}
 \right\}\\
=&0.\end{align*}

So we have proved the equation \ref{Q1'''}. Then the equation \ref{Q1''} is proved.

$\boxempty $

\hspace{1cm}

Combined with the previous argument, we have shown that: the quadratic relation $(Q1)$ is satisfied by any  unitary highest weight $\mathfrak{co}(J)$-module which has the smallest positive Gelfand-Kirillov dimension.

\begin{Rem}\label{value a}

From the computations in our proof, we have

$a=a(J,k)=\left\{
  \begin{array}{ll}
    -(n-1)k+n-\frac{1}{2}, ~(k=0,\frac{1}{2}), &\emph{if}~ J=\Gamma(2n)\\
  -k^{2}-(n-2)|k|+n-1, ~(k~ \mathrm{is~ a~ half ~integer}), &\emph{if} ~J=\Gamma(2n-1)\\
    \frac{n(n+2)}{16}, ~(k=0,1),  &\emph{if}~ J=\mathcal {H}_{n}(\mathbb{R})\\
    \frac{n^{2}-k^{2}}{4},  ~(k=0,1,...), & \emph{if}~ J=\mathcal {H}_{n}(\mathbb{C}) \\
    n^{2}-n-\frac{k}{2}-\frac{k^{2}}{4}, ~ (k=0,1,...), &\emph{if}~ J=\mathcal {H}_{n}(\mathbb{H}) \\
    18, ~(k=0),  & \emph{if}~ J=\mathcal {H}_{3}(\mathbb{O})
  \end{array}
\right.$
i.e., $a=a(J, k)$ only depends on the highest weight $\lambda=\lambda(k)$ given in  Corollary \ref{weight k}.

\hspace{1cm}

\end{Rem}

%


\section{Acknowledgment}
I would like to thank Professor Guowu MENG, for his guidance, encouragement
and many helpful discussions. Also I would like to thank Markus Hunziker, Jingsong Huang,  Jianshu LI and David A. Vogan for several helpful discussions regarding  Gelfand-Kirillov dimension and associated varieties. I would also like to thank the referee for his or her careful reading  of my manuscript and his or her valuable suggestions.


\end{document}